\documentclass[a4paper,11pt,reqno,oneside]{article}

\usepackage{amsmath}
\usepackage{amsthm}
\usepackage{amssymb}
\usepackage{enumerate}
\usepackage{epic,eepic}

\newtheorem*{thma}{Theorem A}
\newtheorem*{thmb}{Theorem B}
\newtheorem{lem}{Lemma}

\newtheorem{prop}[lem]{Proposition}
\newtheorem{coro}[lem]{Corollary}

\theoremstyle{definition}
\newtheorem{ex}[lem]{\emph{Example}}
\newtheorem{rem}[lem]{\emph{Remark}}
\newtheorem*{defi}{\emph{Definition}}

\theoremstyle{remark}
\newtheorem{app-claim}{Claim}

\DeclareMathOperator{\dd}{d}
\newcommand{\Z}{{\mathbb Z}}
\newcommand{\Q}{{\mathbb Q}}
\newcommand{\R}{{\mathbb R}}
\newcommand{\N}{{\mathbb N}}
\newcommand{\Acal}{\mathcal{A}}
\newcommand{\Ccal}{\mathcal{C}}
\newcommand{\coloneq}{\mathrel{\mathop:}=}
\newcommand{\eqcolon}{\mathrel{=\!\!\mathop:}}
\newcommand{\mat}[1]{\boldsymbol{#1}}
\newcommand{\Mcan}{\boldsymbol{M}_{\!\varphi}}

\begin{document}

\title{\textbf{Matrices of 3iet preserving morphisms}}
\author{P.~Ambro\v{z}\quad Z.~Mas\'akov\'a\quad E.~Pelantov\'a \\[5mm]
{\small Doppler Institute 
\& Department of Mathematics}\\
{\small FNSPE, Czech Technical University, Trojanova 13, 120 00 Praha 2, Czech Republic}\\
{\small E-mail: \texttt{petr.ambroz@fjfi.cvut.cz},
\texttt{masakova@km1.fjfi.cvut.cz}},\\{\small \texttt{pelantova@km1.fjfi.cvut.cz}}}
\date{}
\maketitle

\vspace{-2em}

\begin{abstract}
We study matrices of morphisms preserving the family of words
coding 3-interval exchange transformations. It is well known that
matrices of morphisms preserving sturmian words (i.e.\ words
coding 2-interval exchange transformations with the maximal
possible factor complexity) form the monoid
$\{\boldsymbol{M}\in\mathbb{N}^{2\times
2}\;|\;\det\boldsymbol{M}=\pm1\} = \{
\boldsymbol{M}\in\mathbb{N}^{2\times 2}\;|\;
\boldsymbol{M}\boldsymbol{E}\boldsymbol{M}^T =
\pm\boldsymbol{E}\}$, where $\boldsymbol{E} =
(\!\begin{smallmatrix}0&1\\-1&0\end{smallmatrix})$.

We prove that in case of exchange of three intervals, the matrices
preserving words coding these transformations and having the
maximal possible subword complexity belong to the monoid
$\{\boldsymbol{M}\in\mathbb{N}^{3\times 3}\;|\;
\boldsymbol{M}\boldsymbol{E}\boldsymbol{M}^T = \pm\boldsymbol{E},\
 \det\boldsymbol{M}=\pm 1\}$, where $\boldsymbol{E} =
\Big(\!\begin{smallmatrix}0&1&1\\-1&0&1\\-1&-1&0\end{smallmatrix}\Big)$.
\end{abstract}

\section{Introduction}

Sturmian words are the most studied class of infinite aperiodic
words. By their nature, they are defined purely over a binary
alphabet. There exist several equivalent definitions of sturmian
words~\cite{berstel-ijac-12}, which give rise to several different
generalizations of sturmian words over larger alphabets. For
example, the generalization of sturmian words to Arnoux-Rauzy
words of order $r$ uses the characterization of sturmian words by
means of the so-called left and right special
factors~\cite{arnoux-rauzy-bcmf-119}.

Another natural generalization can be derived from the definition
of a sturmian word as an aperiodic word coding a transformation of
exchange of two intervals. The $r$-interval exchange
transformation has been introduced by Katok and
Stepin~\cite{katok-stepin-umn-22}: An exchange $T$ of $r$
intervals is defined by a vector of $r$ lengths and by a
permutation of $r$ letters; the unit interval is then partitioned
according to the vector of lengths, and $T$ interchanges these
intervals according to the given permutation. Rauzy was the first
one to observe that interval exchange transformation can be used
for the generalization of sturmian words.

In contrast to ergodic properties of these transformations, which
were studied by many
authors~\cite{keane-mz-141,rauzy-aa-34,veech-ajm-6,vershik-livshits-asm-9},
combinatorial properties of associated words have been so far
explored only a little. Some results, analogical to the properties
known for sturmian words, have been derived for the most simple
case, namely for 3-interval exchange transformations. Note that
for the exchange of three intervals, the most interesting
permutation is $(321)$ and all the results cited below apply to
transformations with this permutation. Words coding 3-interval
exchange transformation can be periodic or aperiodic, depending on
the choice of parameters. In accordance with the terminology
introduced by~\cite{damanik-zamboni-rmp-15}, infinite words which
code 3-interval exchange transformations and are aperiodic, are
called 3iet words. The factor complexity $\Ccal_u(n)$ of a 3iet
word $u$, i.e., the number of different factors of length $n$
occurring in $u$, is known to satisfy $\Ccal_u(n)\leq 2n+1$ for
all $n\in\N$. Words for which $\Ccal_u(n)= 2n+1$, for all
$n\in\N$, are called non-degenerated (or regular) 3iet words.

In the paper~\cite{ferenczi-holton-zamboni-jam-89}, minimal
sequences coding 3-interval exchange transformations are fully
characterized. The structure of palindromes of these words was
described
in~\cite{damanik-zamboni-rmp-15,balazi-masakova-pelantova-tcs},
whereas the paper~\cite{ferenczi-holton-zamboni-jam-89} deals with
their return words. Here we study morphisms which map the set of
3iet words to itself.


Morphisms preserving sturmian words were completely described by Berstel, Mignosi and
S\'e\'ebold~\cite{berstel-seebold-bbms-1,mignosi-seebold-jtnb-5,seebold-tcs-88}.
Recall that there are two ways to define such a morphism:
\begin{itemize}
\item
  A morphism $\varphi$ over the binary alphabet $\{0,1\}$ is said to be \emph{locally sturmian}
  if there is a sturmian word $u$ such that $\varphi(u)$ is also sturmian.
\item
  A morphism $\varphi$ over the binary alphabet $\{0,1\}$ is said to be \emph{sturmian}
  if $\varphi(u)$ is sturmian for all sturmian words $u$.
\end{itemize}

Berstel, Mignosi and S\'e\'ebold showed that the families of
sturmian and locally sturmian morphisms coincide and that they
form a monoid generated by three morphisms, $\psi_1$, $\psi_2$ and
$\psi_3$, given by
\begin{equation}\label{eq:sturm-generators}
\psi_1:
\begin{aligned}
0 & \mapsto 01 \\ 1 &\mapsto 1
\end{aligned}\,, \quad\quad
\psi_2:
\begin{aligned}
0 & \mapsto 10 \\ 1 &\mapsto 1
\end{aligned}\,, \quad\quad
\psi_3:
\begin{aligned}
0 & \mapsto 1 \\ 1 &\mapsto 0
\end{aligned}\,.
\end{equation}

To each morphism $\varphi$ over a $k$-letter alphabet
$\{a_1,\ldots,a_k\}$ one can assign its incidence matrix
$\Mcan\in\N^{k\times k}$ by putting
\begin{equation}\label{eq:subst-matrix}
(\Mcan)_{ij} = \text{number of letters $a_j$ in the word $\varphi(a_i)$}\,.
\end{equation}

As a simple consequence of the fact that the monoid of sturmian
morphisms is generated by $\psi_1$, $\psi_2$ and $\psi_3$
from~\eqref{eq:sturm-generators}, one has the following fact:
\emph{A matrix $\mat{M}\in\N^{2\times 2}$ is the incidence matrix
of a sturmian morphism if and only if\:\ $\det\mat{M}=\pm1$.} By
an easy calculation we can derive that for matrices of order
$2\times 2$
\[
\det\mat{M}=\pm1 \quad \Longleftrightarrow \quad \mat{M}\mat{E}\mat{M}^T = \pm\mat{E}\,,
\text{ where } \mat{E}=\bigl(\!\begin{smallmatrix} 0 & 1 \\ -1 & 0 \end{smallmatrix}\bigr)\,.
\]
In the theory of Lie groups, one can formulate this claim by
stating that the group $\text{SL}(2,\Z)$ is isomorphic to the
group $\text{Sp}(2,\Z)$, see~\cite{jacobson-liealg}.

The aim of this paper is to derive similar properties for matrices
of morphisms preserving the family of 3iet words, which we call
here 3iet preserving morphisms. We will prove the following
theorems.

\begin{thma}
  Let $\varphi$ be a 3iet preserving morphism and let $\mat{M}$ be its incidence matrix. Then
  \[
  \mat{M}\mat{E}\mat{M}^T = \pm\mat{E}, \text{ where }
  \mat{E}=\Bigl(\!\begin{smallmatrix}0&1&1\\-1&0&1\\-1&-1&0\end{smallmatrix}\Bigr)\,.
  \]
\end{thma}

\begin{thmb}
  Let $\varphi$ be a 3iet preserving morphism and let $\mat{M}$ be its incidence matrix.
  Then one of the following holds \vspace{-0.5em}
  \begin{itemize}
  \item
    $\det\mat{M}=\pm 1$ and $\varphi(u)$ is non-degenerated for every non-degenerated 3iet
    word $u$,\vspace{-0.5em}
  \item
    $\det\mat{M}=0$ and $\varphi(u)$ is degenerated for every 3iet word $u$.
\end{itemize}
\end{thmb}

In the proof of Theorem~A we use the description of matrices of sturmian morphisms
given above, while the main tool employed in the proof of Theorem~B is
the connection between words coding 3-interval exchange transformations and
cut-and-project sets.

Note that we do not address at all the description of the 3iet
preserving morphisms themselves.

\section{Preliminaries}

In this paper we deal with finite and infinite words over a finite
alphabet $\Acal$, whose elements are called letters. The set of
all finite words over $\Acal$ is denoted by $\Acal^*$. This set,
equipped with the concatenation as a binary operation, is a free
monoid having the empty word as its identity. The length of a word
$w=w_1w_2\cdots w_n$ is denoted by $|w|=n$, the number of letters
$a_i$ in the word $w$ is denoted by $|w|_{a_i}$.

\subsection{Infinite words}

The set of two-sided infinite words over an alphabet $\Acal$,
i.e., of two-sided infinite sequences of letters of $\Acal$, is
denoted by $\Acal^\Z$, its elements are words $u=(u_n)_{n\in\Z}$.
Note that in all our considerations we will not identify infinite
words $(u_{n+k})_{n\in\Z}$ and $(u_{n})_{n\in\Z}$, and therefore
we will mark the position corresponding to the index $0$, usually
using $|$ as the delimiter, e.g.\ for $u\in\Acal^\Z$,
\[
u = \cdots u_{-3}u_{-2}u_{-1}|\:u_0u_1u_2\cdots\,.
\]
The words of this form are sometimes called pointed biinfinite words. Naturally, one can
define a metric on the set $\Acal^\Z$.

\begin{defi}
  Let $u=(u_n)_{n\in\Z}$ and $v=(v_n)_{n\in\Z}$ be two biinfinite
  words over $\Acal$. We define the \emph{distance}
  $\dd(u,v)$ between $u$ and $v$ by setting
  \begin{equation}\label{eq:metric}
  \dd(u,v) \coloneq \frac{1}{1+j}\,,
  \end{equation}
  where $j\in\N$ is the minimal index such that either $u_j\neq v_j$ or $u_{-j}\neq v_{-j}$.
\end{defi}
It can be easily verified that the above defined distance
$\dd(u,v)$ is a metric and that the set $\Acal^\Z$ with this
metric is a compact metric space.

We consider also one-sided infinite words $u = (u_n)_{n\in\N}$, either right-sided
$u = u_0u_1u_2\cdots$ or left-sided $u = \cdots u_2u_1u_0$.

The degree of diversity of an infinite word $u$ is expressed by
the complexity function, which counts the number of factors of
length $n$ in the word $u$. Formally, a word $w$ of length $n$ is said
to be a \emph{factor} of a word $u=(u_n)_{n\in\Z}$ if there is an index $i\in\Z$ such
that $w=u_i u_{i+1}\cdots u_{i+n-1}$. The set of all factors of $u$ of length $n$
is denoted by $\mathcal{L}_n(u)$. The \emph{language} $\mathcal{L}(u)$ of an infinite
word $u$ is the set of all its factors, that is,
\[
\mathcal{L}(u) = \bigcup_{n\in\N} \mathcal{L}_n(u)\,.
\]

The \emph{(factor) complexity} $\Ccal_u$ of an infinite word $u$ is the function
$\Ccal_u:\N\rightarrow\N$ defined as
\[
\Ccal_u(n) \coloneq \#\mathcal{L}_n(u)\,.
\]
Clearly, $\Ccal_u(n)$ is a non-decreasing function. Recall that if
there exists $n_0\in\N$ such that $\Ccal_u(n_0)\leq n_0$, then the
word $u$ is eventually periodic (if $u=(u_n)_{n\in\N}$), or
periodic (if $u=(u_n)_{n\in\Z}$), see~\cite{morse-hedlund-ajm-62}.
Hence for an aperiodic word $u$, one has $\Ccal_u(n)\geq n+1$, for
all $n\in\N$.

A one-sided sturmian word $(u_n)_{n\in\N}$ is often defined as an
aperiodic word with complexity $C_u(n)=n+1$, for all $n\in\N$.
However, for biinfinite words, the condition $\Ccal_u(n)\geq n+1$
is not enough for $u$ to be aperiodic. For example, the word
$\cdots111|000\cdots$ has the complexity $\Ccal(n)=n+1$ for all
$n\in\N$. In order to define a biinfinite sturmian word
$(u_n)_{n\in\Z}$ by means of complexity, we need to add another
condition. We introduce the notion of the density of letters,
representing the frequency of occurrence of a given letter in an
infinite word.

The \emph{density of a letter} $a\in\Acal$ in a word
$u\in\Acal^\Z$ is defined as
\[
\rho(a) \coloneq \lim_{n\rightarrow\infty}
\frac{\#\{i\ |\ -n\leq i\leq n,\ u_i=a\}}{2n+1}\,,
\]
if the limit exists.

A biinfinite word $u=(u_n)_{n\in\Z}$ is called sturmian, if
$\Ccal_u(n)=n+1$ for each $n\in\N$ and the densities of letters
are irrational.

Another equivalent definition of sturmian words uses the balance
property. We say that an infinite word $u$ over the alphabet
$\{0,1\}$ is \emph{balanced}, if for every pair of factors
$v,w\in\mathcal{L}_n(u)$ we have $\bigl||v|_0- |w|_0\bigr|\leq 1$. A
one-sided infinite word over the alphabet $\{0,1\}$ is sturmian, if
and only if it is balanced. A biinfinite word over $\{0,1\}$ is
sturmian, if and only if it is balanced and has irrational densities
of letters. For other properties of one-sided and two-sided infinite
sturmian words the reader is referred to~\cite{lothaire2,fogg}.

Unlike the metric space $\Acal^\Z$, the set of all sturmian words
equipped with the same metric~\eqref{eq:metric} is not compact,
however we have the following result.

\begin{lem}\label{lem:limit-of-sturm}
  Let $u\in\{0,1\}^\Z$ be a limit of a sequence of sturmian words $u^{(m)}$.
  Then $u$ is either sturmian or the densities of letters in $u$ are rational.
\end{lem}

\begin{proof}
  Let $w,\widehat{w}\in\mathcal{L}(u)$ be factors of the same length in $u$.
  Since $u = \lim_{m\rightarrow\infty} u^{(m)}$ there exists $m_0\in\N$ such
  that $w,\widehat{w}$ are factors of $u^{(m_0)}$, which is sturmian. Therefore
  $\bigl||w|_0 - |\widehat{w}|_0\bigr|\leq 1$ and $u$ is balanced. If, moreover, the
  densities are irrational, then $u$ is sturmian. The statement follows.
\end{proof}

\subsection{Morphisms and incidence matrices}\label{sec:geometric}

A mapping $\varphi:\Acal^*\rightarrow\Acal^*$ is said to be a
\emph{morphism} over $\Acal$ if
$\varphi(w\widehat{w})=\varphi(w)\varphi(\widehat{w})$ holds for
any pair of finite words $w,\widehat{w}\in\Acal^*$. Obviously, a
morphism is uniquely determined by the images $\varphi(a)$ for all
letters $a\in\Acal$.

The action of a morphism $\varphi$ can be naturally extended to biinfinite words
by the prescription
\[
\varphi(u) = \varphi(\cdots u_{-2}u_{-1}|\:u_0u_1\cdots) \coloneq
\cdots\varphi(u_{-2})\varphi(u_{-1})|\,\varphi(u_{0})\varphi(u_{1})\cdots\,.
\]
The mapping $\varphi : u \mapsto \varphi(u)$ is continuous on $\Acal^\Z$;
a word $u\in\Acal^\Z$ is said to be a \emph{fixed point} of $\varphi$ if $\varphi(u)=u$.

Recall that the incidence matrix of a morphism $\varphi$ over the
alphabet $\Acal$ is defined by~\eqref{eq:subst-matrix}.
A morphism $\varphi$ is called primitive if there exist an integer $k$ such that
the matrix $\mat{M}_{\varphi}^k$ is positive.

Morphisms over $\Acal$ form a monoid, whose neutral element is the
identity morphism. Let $\varphi$ and $\psi$ be morphisms over
$\Acal$, then the matrix of their composition, that is, of the
morphism $u\mapsto(\varphi\circ\psi)(u) =
\varphi\bigl(\psi(u)\bigr)$ is obtained by
\begin{equation}\label{eq:matrix-compose}
\mat{M}_{\varphi\circ\psi} = \mat{M}_{\psi}\mat{M}_{\varphi}\,.
\end{equation}

Let us now explain the importance of the incidence matrix of a
morphism $\varphi$ for the combinatorial properties of infinite
words on which the morphism $\varphi$ acts. Assume that an
infinite word $u$ over the alphabet $\Acal=\{a_1,\dots, a_k\}$ has
well defined densities of letters, given by the vector
\[
\vec{\rho}_u = \bigl(\rho(a_1),\dots,\rho(a_k)\bigr)\,.
\]
It is easy to see that the densities of letters in the infinite
word $\varphi(u)$ are also well defined and it
holds that
\begin{equation}\label{eq:density-phi-u}
  \vec{\rho}_{\varphi(u)} = \frac{\vec{\rho}_u\Mcan}
  {\vec{\rho}_u\Mcan
  \Bigl(\begin{smallmatrix}1\\[-2mm] \vdots\\1\end{smallmatrix}\Bigr)}\,,
\end{equation}
where $\Mcan$ is the incidence matrix of $\varphi$.

Assume now that the infinite word $u$ is a fixed point of a
morphism $\varphi$. Then from~\eqref{eq:density-phi-u}, we obtain
that the vector of densities $\vec{\rho}_u$ is a left eigenvector
of the incidence $\Mcan$, i.e.,
$\vec{\rho}_u\Mcan=\Lambda\vec{\rho}_u$. Since $\Mcan$ is a
non-negative integral matrix, we can use the Perron-Frobenius
Theorem stating that $\Lambda$ is the dominant eigenvalue of
$\Mcan$. Moreover, all eigenvalues of $\Mcan$ are algebraic
integers.

The right eigenvector of the incidence matrix corresponding to the
dominant eigenvalue has also a nice interpretation. It plays an
important role for the \emph{geometric representation} of a fixed
point of a morphism. Let $u$ be a fixed point of a morphism
$\varphi$ over a $k$-letter alphabet $\{a_1,\ldots,a_k\}$ and let
$\Mcan$ have a positive right eigenvector $\vec{x}$. The infinite
word $u$ can be geometrically represented by a self-similar set
$\Sigma$ as follows.

Let us denote by $x_1,x_2,\ldots,x_k$ the positive components of
$\vec{x}$, and let $\Lambda$ be the corresponding eigenvalue,
i.e., $\Mcan\vec{x}=\Lambda\vec{x}$. Since $\Mcan$ is non-negative
and $\vec{x}$ is positive, the eigenvalue $\Lambda$ is equal to
the spectral radius of the matrix $\Mcan$. Moreover, $\Mcan$ being
an integral matrix implies $\Lambda\geq 1$.

For a biinfinite word $u = \cdots
u_{-3}u_{-2}u_{-1}|u_0u_1u_2\cdots$ we denote
\begin{multline*}
  \qquad \Sigma = \Bigl\{ \sum_{i=1}^k |w|_{a_i}x_i\ \Big|\
  w \text{ is an arbitrary prefix of } u_0u_1u_2\cdots \Bigr\} \\
  \cup \Bigl\{ -\sum_{i=1}^k |w|_{a_i}x_i\ \Big|\
  w \text{ is an arbitrary suffix of } \cdots u_{-3}u_{-2}u_{-1} \Bigr\}\,. \qquad
\end{multline*}
The set $\Sigma$ can be equivalently defined as
\[
\Sigma = \{t_n\ |\ n\in\Z\}\,,\quad \text{where} \quad t_0=0 \
\hbox{ and } \ t_{n+1} - t_n = x_i\ \Leftrightarrow\ u_n=a_i\,.
\]
Since $u$ is a fixed point of a morphism, the construction of
$\Sigma$ implies that $\Lambda\Sigma \subset \Sigma$.
A set having this property is called self-similar.
\begin{figure}[!ht]
\begin{center}
\begin{picture}(370,155)
\put(4,28){\line(1,0){364}}
\put(4,100){\line(1,0){364}}
\put(7,23){\line(0,1){82}}
    \put(7,28){\circle*{5}}
    \put(23,20){\makebox(0,0){$1$}} \put(23,107){\makebox(0,0){$1$}}
\put(40,23){\line(0,1){82}}
    \put(56,20){\makebox(0,0){$1$}} \put(56,107){\makebox(0,0){$1$}}
\put(72,23){\line(0,1){82}}
    \put(82,20){\makebox(0,0){$0$}} \put(82,107){\makebox(0,0){$0$}}
\put(92,20){\line(0,1){85}}
    \put(92,28){\circle*{5}} \put(92,100){\circle*{5}}
    \put(108,20){\makebox(0,0){$1$}} \put(112,106){\makebox(0,0){$1$}}
\put(124,23){\line(0,1){82}}
    \put(124,100){\circle*{5}}
    \put(134,20){\makebox(0,0){$0$}} \put(136,107){\makebox(0,0){$0$}}
\put(144,20){\line(0,1){125}}
    \put(144,100){\circle*{5}} \put(144,144.45){\circle*{5}} \put(144,28){\circle*{5}}
    \put(160,20){\makebox(0,0){$1$}} \put(160,107){\makebox(0,0){$1$}}
\put(176,23){\line(0,1){82}}
    \put(176,100){\circle*{5}}
    \put(192,20){\makebox(0,0){$1$}} \put(189,106){\makebox(0,0){$1$}}
\put(208,23){\line(0,1){82}}
    \put(208,100){\circle*{5}}
    \put(218,20){\makebox(0,0){$0$}} \put(218,107){\makebox(0,0){$0$}}
\put(228,20){\line(0,1){85}}
    \put(228,100){\circle*{5}} \put(228,28){\circle*{5}}
    \put(244,20){\makebox(0,0){$1$}} \put(244,107){\makebox(0,0){$1$}}
\put(260,23){\line(0,1){82}}
    \put(276,20){\makebox(0,0){$1$}} \put(276,107){\makebox(0,0){$1$}}
\put(292,23){\line(0,1){82}}
    \put(302,20){\makebox(0,0){$0$}} \put(302,107){\makebox(0,0){$0$}}
\put(312,23){\line(0,1){82}}
    \put(312,28){\circle*{5}}
    \put(328,20){\makebox(0,0){$1$}} \put(328,107){\makebox(0,0){$1$}}
\put(344,23){\line(0,1){82}}
    \put(354,20){\makebox(0,0){$0$}} \put(354,107){\makebox(0,0){$0$}}
\put(364,23){\line(0,1){82}}
    \put(364,28){\circle*{5}}
\put(160,88){\makebox(0,0){$\underbrace{\hspace*{30pt}}_{\hbox{\footnotesize $\ell(1)$}}$}}
\put(134,89){\makebox(0,0){$\underbrace{}_{\hbox{\footnotesize $\,\ell(0)$}}$}}
\put(186,5){\makebox(0,0){$\underbrace{\hspace*{83pt}}_{\hbox{$\lambda\ell(1)$}}$}}
\put(118,5){\makebox(0,0){$\underbrace{\hspace*{51pt}}_{\hbox{$\lambda \ell(0)$}}$}}
\qbezier(144,144.45)(92,100)(7,28)
\qbezier(144,144.45)(124,100)(92,28)
\qbezier(144,144.45)(176,100)(228,28)
\qbezier(144,144.45)(208,100)(312,28)
\qbezier(144,144.45)(228,100)(364,28)
\end{picture}
\caption{Action of the morphism $0\mapsto10$, $1\mapsto110$ on the
geometric representation of its fixed point
$u=\lim_{n\rightarrow\infty}\varphi^n(0)|\varphi^n(1)$.}
\label{fig:morphism}
\end{center}
\end{figure}
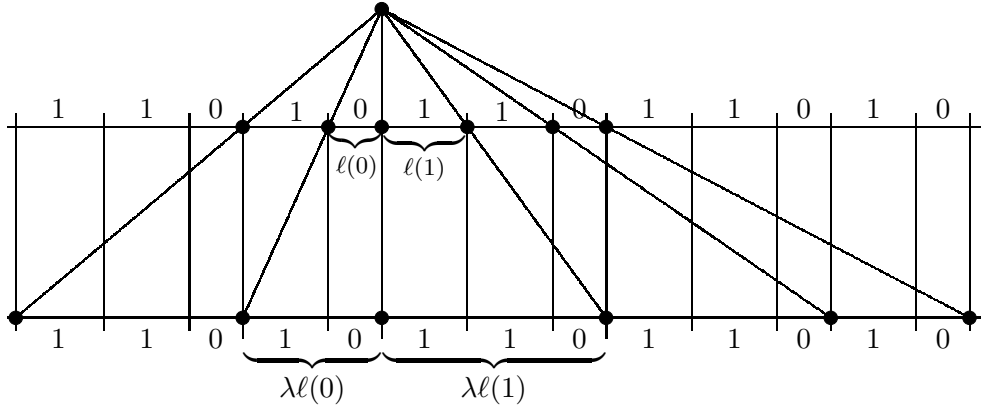

Moreover, if $u_n=a_i$ then the number of points of the set
$\Sigma$ belonging to $(\Lambda t_n,\Lambda t_{n+1}]$ is equal to
the length of $\varphi(a_i)$. Formally, we have
\begin{equation}\label{eq:points-in-tntn1}
  \# \bigl((\Lambda t_n,\Lambda t_{n+1}] \cap \Sigma\bigr) = |\varphi(a_i)|\,.
\end{equation}

In Figure~\ref{fig:morphism}, one can see the geometric
representation of the fixed point of the morphism $0\mapsto10$,
$1\mapsto110$. The matrix of this morphism, $\mat{M} =
\big(\begin{smallmatrix}1&1\\1&2\end{smallmatrix}\big)$, has the
dominant eigenvalue $\Lambda=\tau^2$, where
$\tau=\tfrac{1+\sqrt{5}}{2}$ is the golden ratio. The
corresponding right eigenvector of $M$ is
$\big(\begin{smallmatrix}1\\\tau\end{smallmatrix}\big)$. Hence the
lengths assigned to letters $0$ and $1$ are $\ell(0)=1$ and
$\ell(1)=\tau$, respectively.

\section{Interval exchange words}

Before we define infinite words coding a 3-interval exchange
transformation, we will show the definition of sturmian words
using a 2-interval exchange transformation. It is well known (see
e.g.~\cite{morse-hedlund-ajm-62,lothaire2}) that every sturmian
word $u=(u_n)_{n\in\Z}$ over the alphabet $\{0,1\}$ satisfies
\begin{equation}\label{eq:dolni-mechanicke}
u_n = \lfloor (n+1)\alpha + x_0\rfloor - \lfloor n\alpha+x_0\rfloor
\qquad \text{for all $n\in\Z$,}
\end{equation}
 or
\begin{equation}\label{eq:horni-mechanicke}
u_n = \lceil (n+1)\alpha + x_0\rceil - \lceil n\alpha+x_0\rceil
\qquad \text{for all $n\in\Z$,}
\end{equation}
where $\alpha\in(0,1)$ is an irrational number called the slope,
and $x_0\in[0,1)$ is called the intercept of $u$. In the former
case, $(u_n)_{n\in\Z}$ is the so-called upper mechanical word, in
the latter case the lower mechanical word, with slope $\alpha$ and
intercept $x_0$.

If $(u_n)_{n\in\Z}$ is of the form~\eqref{eq:dolni-mechanicke} then, obviously,
\begin{equation}\label{eq:un-dolni-mech}
u_n = \begin{cases}
0 & \text{if $\{n\alpha+x_0\}\in[0,1-\alpha)\,,$} \\
1 & \text{if $\{n\alpha+x_0\}\in[1-\alpha,1)\,,$}
\end{cases}
\end{equation}
where $\{x\}$ denotes the fractional part of $x$, i.e.,
$\{x\}=x-\lfloor x\rfloor$. We can define a transformation
$T:[0,1)\rightarrow[0,1)$ by the prescription
\begin{equation}\label{eq:2iet}
T(x) = \begin{cases}
x + \alpha & \text{if $\{n\alpha+x_0\}\in [0,1-\alpha)\eqcolon I_0\,,$} \\
x + \alpha - 1 & \text{if $\{n\alpha+x_0\}\in [1-\alpha,1)\eqcolon I_1\,,$}
\end{cases}
\end{equation}
which satisfies $T(x) = \{x+\alpha\}$. It follows easily that the $n$-th
iteration of $T$ is given as
\begin{equation}\label{eq:T-n}
T^n(x) = \{x+n\alpha\} \qquad \text{for all $n\in\Z$}.
\end{equation}
Putting~\eqref{eq:un-dolni-mech} and~\eqref{eq:T-n} together, we see that
a sturmian word $(u_n)_{n\in\Z}$ can be defined using the transformation
$T$ by
\[
u_n = \begin{cases}
0 & \text{if $T^n(x_0)\in I_0$,} \\
1 & \text{if $T^n(x_0)\in I_1$.}
\end{cases}
\]
Hence a sturmian word is given by iterations of the intercept $x_0$
under the mapping $T$, that is, by the orbit of $x_0$ under $T$.

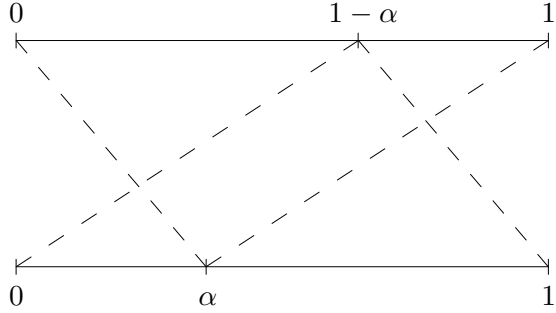
\begin{figure}[!ht]
\centering
\setlength{\unitlength}{1cm}
\begin{picture}(9,4.5)(0,0.5)
\drawline(1,4)(8,4)
\drawline(1,3.9)(1,4.1) \drawline(5.5,3.9)(5.5,4.1) \drawline(8,3.9)(8,4.1)
\drawline(1,1)(8,1)
\drawline(1,0.9)(1,1.1) \drawline(3.5,0.9)(3.5,1.1) \drawline(8,0.9)(8,1.1)
\dashline{0.25}(1,4)(3.5,1) \dashline{0.25}(5.5,4)(8,1)
\dashline{0.25}(5.5,4)(1,1) \dashline{0.25}(8,4)(3.5,1)
\put(0.91,0.5){0} \put(3.4,0.5){$\alpha$} \put(7.91,0.5){1}
\put(0.91,4.25){0} \put(5.11,4.25){$1-\alpha$} \put(7.91,4.25){1}
\end{picture}
\caption{Graph of a 2-interval exchange transformation.}\label{fig:2iet}
\end{figure}

The action of the mapping $T$ from~\eqref{eq:2iet} is illustrated
on Figure~\ref{fig:2iet}.

We see that $T$ is in fact an exchange of two intervals $I_0=[0,1-\alpha)$ and
$I_1=[1-\alpha,1)$. It is therefore called a 2-interval exchange transformation.

Let us mention that if $(u_n)_{n\in\Z}$ is an upper mechanical word, the corresponding
2-interval exchange transformation is given by
$T:(0,1]\mapsto(0,1]$, with $I_0=(0,1-\alpha]$ and $I_1=(1-\alpha,1]$.
Note also that it was not necessary that $T$ was acting on a unit interval.
We could choose an arbitrary interval divided into two parts, ratio of whose
lengths would be irrational.

\medskip
Analogically to the case of exchange of two intervals, we can
define a 3-interval exchange transformation.

\begin{defi}
  Let $\alpha,\beta,\gamma$ be three positive real numbers. Denote
  $$
  \begin{array}{rcl}
  I_A&\coloneq&[0,\alpha)\\
  I_B&\coloneq& [\alpha,\alpha+\beta)\\
  I_C&\coloneq& [\alpha+\beta,\alpha+\beta+\gamma)
  \end{array}
  \quad\hbox{or }\quad
  \begin{array}{rcl}
  I_A&\coloneq& (0,\alpha]\\
  I_B&\coloneq& (\alpha,\alpha+\beta]\\
  I_C&\coloneq& (\alpha+\beta,\alpha+\beta+\gamma]
  \end{array}
  $$
  respectively, and $I:=I_A\cup I_B\cup I_C$.
  A mapping
  $T:I\rightarrow I$, given by
  \begin{equation}\label{eq:3iet}
  T(x) = \begin{cases}
  x + \beta+\gamma & \text{if $x\in I_A$,} \\
  x -\alpha+\gamma & \text{if $x\in I_B$,} \\
  x -\alpha-\beta & \text{if $x\in I_C$,}
  \end{cases}
  \end{equation}
  is called a \emph{3-interval exchange transformation} (3iet)\footnote{Note that the
  above defined mapping $T$ should be more precisely called 3-interval
  exchange with the permutation $(321)$, since the initial arrangement
  of intervals $I_A<I_B<I_C$ is changed to $T(C)<T(B)<T(A)$. Indeed,
  one can define also 3iet with a different permutation of intervals,
  e.g.\ (231). The corresponding 3iet word has the property that by
  changing all the letters $C$ into $B$ one obtains a sturmian word
  over the alphabet $\{A,B\}$. We will not consider such words.}
  with parameters $\alpha,\beta,\gamma$.
\end{defi}

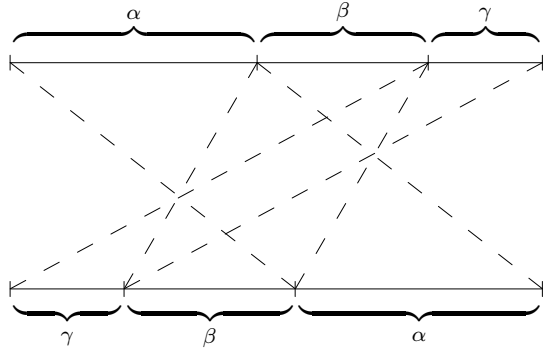
\begin{figure}[ht!]
\centering
\setlength{\unitlength}{1cm}
\begin{picture}(9,4.6)(0,0.4)
\drawline(1,4)(8,4)
\drawline(1,3.9)(1,4.1) \drawline(4.25,3.9)(4.25,4.1)
\drawline(6.5,3.9)(6.5,4.1) \drawline(8,3.9)(8,4.1)
\drawline(1,1)(8,1)
\drawline(1,0.9)(1,1.1) \drawline(2.5,0.9)(2.5,1.1)
\drawline(4.75,0.9)(4.75,1.1) \drawline(8,0.9)(8,1.1)
\dashline{0.25}(1,4)(4.75,1) \dashline{0.25}(4.25,4)(8,1)
\dashline{0.25}(4.25,4)(2.5,1) \dashline{0.25}(6.5,4)(4.75,1)
\dashline{0.25}(6.5,4)(1,1) \dashline{0.25}(8,4)(2.5,1)
\put(1.05,4.25){$\overbrace{\hspace{3.15cm}}^{\alpha}$}
\put(4.3,4.25){$\overbrace{\hspace{2.15cm}}^{\beta}$}
\put(6.55,4.25){$\overbrace{\hspace{1.4cm}}^{\gamma}$}
\put(1.05,0.8){$\underbrace{\hspace{1.4cm}}_{\gamma}$}
\put(2.55,0.8){$\underbrace{\hspace{2.15cm}}_{\beta}$}
\put(4.8,0.8){$\underbrace{\hspace{3.15cm}}_{\alpha}$}
\end{picture}
\caption{Graph of a 3-interval exchange transformation.}\label{fig:3iet}
\end{figure}

The graph of a 3-interval exchange is on Figure~\ref{fig:3iet}.

With a 3-interval exchange transformation $T$, one can naturally
associate a ternary biinfinite word $u_T(x_0)=(u_n)_{n\in\Z}$,
which codes the orbit of a point $x_0$ from the domain of $T$, as
\begin{equation}\label{eq:u_T}
u_n = \begin{cases}
A & \text{if $\ T^n(x_0)\in I_A$,} \\
B & \text{if $\ T^n(x_0)\in I_B$,} \\
C & \text{if $\ T^n(x_0)\in I_C$.}
\end{cases}
\end{equation}

Similarly as in the case of a 2-interval exchange transformation,
the infinite word coding a 3iet can be periodic or aperiodic,
according to the choice of parameters $\alpha,\beta,\gamma$. We
will focus only on aperiodic words.

\begin{defi}
  An aperiodic\footnote{A biinfinite word $(u_n)_{n\in\Z}$ is called aperiodic
  if neither $u_0u_1u_2\cdots$ nor $\cdots u_{-3}u_{-2}u_{-1}$ is eventually
  periodic.} word $u_T(x_0)$ coding the orbit of the point $x_0$
  under the 3iet $T$ defined above is called a \emph{3iet word} with parameters
  $\alpha,\beta,\gamma$ and $x_0$.
\end{defi}

The following lemma shows a close relation between words coding
3-interval exchange and 2-interval exchange transformations.
\begin{lem}\label{lem:delta-of-3iet}
  Let $u = (u_n)_{n\in\Z}$ be a word coding 3-interval exchange transformation and
  let $\sigma:\{A,B,C\}^*\rightarrow\{0,1\}^*$ be a morphism given by
  \begin{equation}\label{eq:delta}
  A \mapsto 0\,, \qquad
  B \mapsto 01\,, \qquad
  C \mapsto 1\,.
  \end{equation}
  Then $\sigma(u)$ codes a 2-interval exchange transformation.
\end{lem}

\begin{proof}
  Let $u$ be the coding of $x_0$ under the 3-interval exchange transformation $T$ with
  intervals $[0,\alpha)$, $[\alpha,\alpha+\beta)$
  and $[\alpha+\beta,\alpha+\beta+\gamma)$.

  Let $S$ be the 2-interval exchange transformation of the intervals
  $I_0=[0,\alpha+\beta)$ and $I_1=[\alpha+\beta,\alpha+2\beta+\gamma)$, i.e.,
  \[
  S(x) = \begin{cases}
  x + \beta  + \gamma & \text{if $x\in I_0$,} \\
  x - \alpha - \beta & \text{if $x\in I_1$.}
  \end{cases}
  \]
  One can easily see that
  \begin{align*}
  x & \in[0,\alpha) && \Rightarrow \qquad x\in I_0 \text{ and } T(x)=S(x)\,, \\
  x & \in[\alpha,\alpha+\beta) && \Rightarrow \qquad x\in I_0,\ S(x)\in I_1
  \text{ and } S^2(x) = T(x)\,, \\
  x & \in[\alpha+\beta,\alpha+\beta+\gamma) && \Rightarrow \qquad x\in I_1,
  \text{ and } S(x) = T(x)\,.
  \end{align*}
  This proves that $\sigma(u)$ is the coding of $x_0$ under $S$.
\end{proof}

\section{Periodic and aperiodic words coding 3iet}

In order to clarify the relation between the parameters of a 3iet
and the complexity of the corresponding infinite words, we recast
the definition of these words in a new formalism. We show that
every 3iet word codes distances in a discrete set arising as a
projection of points of the lattice $\Z^2$. This construction is
known as the cut-and-project method.

Let $\varepsilon,\eta$ be real numbers, $\varepsilon\neq-\eta$.
Every point $(a,b)\in\Z^2$ can be written in the form
\[
(a,b) = (a+b\eta)\vec{x}_1 + (a-b\varepsilon)\vec{x}_2\,,
\]
where
\[
\vec{x}_1 = \frac{1}{\varepsilon+\eta}(\varepsilon,1)
\qquad\text{and}\qquad \vec{x}_2 =
\frac{1}{\varepsilon+\eta}(\eta,-1)\,.
\]

Let $V_1$ and $V_2$ denote the lines in $\R^2$ spanned by
$\vec{x}_1$ and $\vec{x}_2$, respectively. Then
$(a+b\eta)\vec{x}_1$ is the projection of the lattice point
$(a,b)$ on $V_1$, whereas $(a-b\varepsilon)\vec{x}_1$ is its projection
on $V_2$. Let $\Omega$ be a bounded interval. Then the set
\begin{equation}\label{eq:Sigma}
\Sigma_{\varepsilon,\eta}(\Omega) \coloneq \{ a+b\eta\ |\
a,b\in\Z, \ a-b\varepsilon\in\Omega\}
\end{equation}
is called the \emph{Cut-and-project (C\&P) set} with parameters
$\varepsilon,\eta,\Omega$. Thus C\&P sets arise by projection on
the line $V_1$ of points of $\Z^2$ having their second projection
in a chosen segment on $V_2$.

\begin{prop}\label{prop:cap}
  Let $\alpha, \beta, \gamma$ be positive real numbers, and
  let $T : [0,\alpha+\beta+\gamma) \mapsto
  [0,\alpha+\beta+\gamma) $ be a 3iet defined by~\eqref{eq:3iet}.
  Let $x_0 \in [0,\alpha+\beta+\gamma)$ and let $u_T(x_0) =
  (u_n)_{n\in \mathbb{Z}}$ be the biinfinite word given by~\eqref{eq:u_T}.
  Put
  \begin{equation}\label{eq:prevod}
    \varepsilon \coloneq
    \frac{\beta+\gamma}{\alpha+2\beta+\gamma}\,,\quad l\coloneq
    \frac{\alpha+\beta+\gamma}{\alpha+2\beta+\gamma}\,,\quad c
    \coloneq \frac{x_0}{\alpha+2\beta+\gamma}\, \quad\hbox{and}\quad
    \Omega = (c-l,c]\,,
  \end{equation}
  and choose arbitrary $\eta >0 $. Then the C\&P set
  $\Sigma_{\varepsilon,\eta}(\Omega)$  is a discrete set with the following properties:
  \begin{enumerate}
  \item
    $0 \in \Sigma_{\varepsilon,\eta}(\Omega)$;
  \item
    the distances between adjacent elements of $\Sigma_{\varepsilon,\eta}(\Omega)$
    take values $\mu_A = \eta$, $\mu_B= 1+2\eta$, and $\mu_C = 1+\eta $;
  \item
    the ordering of the distances with respect to the origin is coded
    by the word $u_T(x_0)$;
  \item
    $\Sigma_{\varepsilon,\eta}(\Omega) = \bigl\{ \lfloor
    c+n\varepsilon\rfloor+n\eta\ |\ n\in\Z,\
    \{c+n\varepsilon\}\in[0,l) \bigr\}$.
  \end{enumerate}
\end{prop}

\begin{proof}
The parameters $\varepsilon$, $l$, and $c$ satisfy clearly
\begin{equation}\label{eq:param-cap}
\varepsilon \in (0,1)\,,\quad \max\{\varepsilon,1-\varepsilon\} <
l \leq 1\,,\quad 0 \in (c-l,c]\,.
\end{equation}
The condition in~\eqref{eq:Sigma} determining whether a given point $a+b\eta$ belongs to the
C\&P set $\Sigma_{\varepsilon,\eta}(\Omega)$ can be rewritten
\[
a-b\varepsilon\in\Omega\quad \Leftrightarrow\quad c+b\varepsilon-l< a\leq
c+b\varepsilon\quad \Leftrightarrow\quad a=\lfloor
c+b\varepsilon\rfloor\text{ and }\{c+b\varepsilon\}\in[0,l)\,.
\]
Therefore, the C\&P set $\Sigma_{\varepsilon,\eta}(\Omega)$ can be
expressed as
\begin{equation}\label{eq:CnP-cl}
\Sigma_{\varepsilon,\eta}(\Omega) = \bigl\{ \lfloor
c+n\varepsilon\rfloor+n\eta\ |\ n\in\Z,\
\{c+n\varepsilon\}\in[0,l) \bigr\}\,.
\end{equation}
 Let us denote
$y_n\coloneq\lfloor c+n\varepsilon\rfloor+n\eta$ and
$y_n^*\coloneq\{c+n\varepsilon\}$. From the choice of the
parameter $\varepsilon$ and $\eta$ we can derive that the sequence
$(y_n)_{n\in\Z}$ is strictly increasing. Since
$\Sigma_{\varepsilon,\eta}(\Omega) \subset \{ y_n\:|\:n\in
\mathbb{Z}\}$, to every element $y\in
\Sigma_{\varepsilon,\eta}(\Omega)$ corresponds a point $y^*$. We
show that the distance of $y$ and its right neighbour depends on
the position of $y^*$ in the interval $[0,l)$. Moreover, if $z$ is
the right neighbour of $y$ in $\Sigma_{\varepsilon,\eta}(\Omega)$,
then $z^*=\widetilde{T}(y^*)$, where
$\widetilde{T}:[0,l)\rightarrow[0,l)$ is a 3iet given by the
prescription
\begin{equation}\label{eq:T-CaP}
\widetilde{T}(x) =
\begin{cases}
 x + \varepsilon & \text{if }x\in[0,l-\varepsilon) \,, \\
 x + 2\varepsilon-1 & \text{if }x\in[l-\varepsilon,1-\varepsilon) \,, \\
 x + \varepsilon-1 & \text{if }x\in[1-\varepsilon,l) \,.
\end{cases}
\end{equation}

Let us determine the right neighbour of a point
$y\in\Sigma_{\varepsilon,\eta}(c-l,c]$. Let $y=y_n$, $n\in\Z$,
i.e., $y_n^*=\{c+n\varepsilon\}\in[0,l)$. We discuss three
separate cases, all the time using the fact that
$\max\{\varepsilon,1-\varepsilon\} < l \leq 1$.
\begin{enumerate}[i)]
\item
  if $y_n^*\in[0,l-\varepsilon)$ then $y_{n+1}^*=\{c+(n+1)\varepsilon\} =
  y_n^* + \varepsilon\in[0,l)$ and
  $\lfloor c+n\varepsilon\rfloor = \lfloor c+(n+1)\varepsilon\rfloor$. Hence the
  distance between $y_n$ and its right neighbour is $y_{n+1}-y_n=\eta$.
\item
  if $y_n^*\in[l-\varepsilon,1-\varepsilon)$ then
  $y_{n+1}^*=\{c+(n+1)\varepsilon\} = y_n^* + \varepsilon \in[l,1)$,
  hence $y_{n+1}$ does not belong to the set $\Sigma_{\varepsilon,\eta}(c-l,c]$.
  However, $y_{n+2}^*=\{c+(n+2)\varepsilon\} = y_n^* + 2\varepsilon - 1\in[0,l)$
  and $\lfloor c+(n+2)\varepsilon\rfloor = 1+\lfloor c+n\varepsilon\rfloor$.
  Therefore the right neighbour of $y_n$ is $y_{n+2}$ and we have
  $y_{n+2}-y_n = 1+2\eta$.
\item
  if $y_n^*\in[1-\varepsilon,l)$ then $y_{n+1}^*=\{c+(n+1)\varepsilon\} =
  y_n^* + \varepsilon-1\in[0,l)$, $y_{n+1}$ is the right neighbour of $y_n$
  and $y_{n+1}-y_n = 1 + \eta$.
\end{enumerate}

As $y_0 = 0 \in \Sigma_{\varepsilon,\eta}(c-l,c]$ and $y_0^* =
\{c\}=c$, the distances between consecutive elements of the C\&P set
$\Sigma_{\varepsilon,\eta}(c-l,c]$ are coded by the infinite word
$u_{\widetilde{T}}(c)$. It is easy to see that with our choice of $l$,
$\varepsilon$, and $c$, the lengths of the partial intervals in the definition
of the 3iet $\widetilde{T}$ and the starting point $c$ are only
$(\alpha+2\beta+\gamma)$-multiples of the partial intervals of the
3iet $T$ and its starting point $x_0$, ($\widetilde{T}$ and $T$ are homothetic 3iets).
Therefore $u_{\widetilde{T}}(c) = u_T(x_0)$.
\end{proof}

Let us mention that a 3iet $T$ with the domain
$(0,\alpha+\beta+\gamma]$ corresponds also to a C\&P set with
parameters similar to~\eqref{eq:prevod}.

It is known that a word coding an $r$-interval exchange
transformation with arbitrary permutation of intervals has
complexity ${\cal C}(n)\leq (r-1)n+1$ for all $n\in\N$,
see~\cite{keane-mz-141}. It is useful to distinguish the words with full
complexity and the others.

\begin{defi}
  A 3iet word is called \emph{non-degenerated}, if ${\cal C}(n)= 2n+1$ for
  all $n\in\N$. Otherwise it is called \emph{degenerated}.
\end{defi}

The following proposition allows one to classify the words coding
3iet according to the parameters to periodic, 3iet degenerate, and
3iet non-degenerate infinite words.

\begin{prop}\label{prop:dege}
  Let $T$ be a 3iet transformation of the interval $I$ with
  parameters $\alpha$, $\beta$, $\gamma$, and let $x_0\in I$.
  \begin{itemize}
  \item
    The infinite word $u_T(x_0)$ defined by~\eqref{eq:u_T} is aperiodic
    if and only if
    $$
      \alpha+\beta \text{ and } \beta+\gamma \text{ are linearly independent over $\Q$.}
    $$
  \item
    If the word $u_T(x_0)$ is aperiodic then it is degenerated if and only if
    $$
      \alpha+\beta+\gamma \in (\alpha+\beta)\Z + (\beta+\gamma)\Z\,.
    $$
  \end{itemize}
\end{prop}

\begin{proof}
The formula~\eqref{eq:CnP-cl} for the C\&P set
$\Sigma_{\varepsilon,\eta}(c-l,c]$ implies easily that if
$\varepsilon$ is rational, then the set
$\Sigma_{\varepsilon,\eta}(\Omega)$ is periodic, i.e., the orbit
of every point under the 3iet $\widetilde{T}$ is periodic. On the
other hand, if $\varepsilon$ is irrational, the sequence
$\{c+n\varepsilon\}$ is uniformly distributed, and thus also the
orbit of every point under $\widetilde{T}$ is dense in $[0,l)$.
The relation~\eqref{eq:prevod} between the parameters
$\varepsilon$ and $\alpha, \beta, \gamma$ implies the statement
about periodicity of $u_T(x_0)$.

The complexity of an infinite word coding a C\&P set with
irrational parameters $\varepsilon,\eta$ has been described
in~\cite{gmp-jtnb-15}. It is shown that such a word has the
complexity ${\cal C}(n)=2n+1$ for all $n$ if and only if the
length $l$ of the interval $\Omega$ from~\eqref{eq:Sigma}
satisfies $l\notin\Z+\Z\varepsilon$. The
relation~\eqref{eq:prevod} implies the necessary and sufficient
condition for the degeneracy of the corresponding infinite word.
\end{proof}

We will use the following reformulation of the above statements.

\begin{coro}\label{coro:period}
The infinite word $u_T(x_0)$, defined by~\eqref{eq:u_T}, with
parameters $\alpha,\beta,\gamma>0$ is
\begin{itemize}
\item
  periodic if there exist $K,L\in\Z$, \ $K,L\neq 0$ such that
  \begin{equation}\label{eq:podminka-periodic}
  (\alpha,\beta,\gamma)\left(\begin{smallmatrix}K\\K+L\\L\end{smallmatrix}\right)= 0\,,
  \end{equation}
\item
  aperiodic degenerate if there exist unique $K,L\in\Z$ such that
  \begin{equation}\label{eq:podminka-degenerate}
  (\alpha,\beta,\gamma)\left(\begin{smallmatrix}1\\1\\1\end{smallmatrix}\right)=
  (\alpha,\beta,\gamma)\left(\begin{smallmatrix}K\\K+L\\L\end{smallmatrix}\right)\,.
  \end{equation}
\end{itemize}
\end{coro}

Note that the sequence $\{c+n\varepsilon\}$ being uniformly
distributed for $\varepsilon$ irrational implies not only the
aperiodicity of the infinite word, but also that the densities of
letters are well defined.

\begin{coro}\label{coro:hustoty}
All letters in a 3iet word $u$ with parameters
$\alpha,\beta,\gamma$ have a well defined density and the vector
of densities of $u$, denoted by $\vec{\rho}_u \coloneq
\big(\rho(A),\rho(B),\rho(C)\big)$, is proportional to the vector
$(\alpha,\beta,\gamma)$.
\end{coro}

For the transformation $T$ of exchange of $r$ intervals,  it is
generally difficult to describe the conditions under which the
corresponding dynamical system is minimal, i.e., under which
condition the orbit $\{T^n(x_0)\ |\ n\in \mathbb{Z}\}$ of any
point $x_0$ is dense in the domain of $T$. Keane provides
in~\cite{keane-mz-141} two sufficient conditions for the
minimality of $T$: one of them is the linear independence of
parameters $\alpha$, $\beta$ and $\gamma$ over ${\mathbb Q}$;
second, weaker condition is that the orbits of all discontinuity
points of $T$ are disjoint. This condition is called i.d.o.c.
In~\cite{ferenczi-holton-zamboni-jam-89} it is shown that the
parameters $\alpha,\beta,\gamma$ fulfill i.d.o.c. if and only if
they satisfy neither~(\ref{eq:podminka-periodic})
nor~(\ref{eq:podminka-degenerate}). Nevertheless, even the weaker
condition i.d.o.c. is only sufficient, but not necessary for the
minimality of the dynamical system of $T$. The geometric
representation of 3iet $T$ using a cut-and-project set allows us
to provide a simple characterization of minimal dynamical systems
among 3iet.

\begin{coro}
The dynamical system given by a 3-interval exchange transformation
$T$ with parameters $\alpha, \beta,\gamma$ is minimal if and only if the numbers
$\alpha + \beta$  and $\beta + \gamma$ are linearly independent over
$\mathbb{Q}$.
\end{coro}

\begin{rem} It can be shown,
(see~\cite{arnoux-rauzy-bcmf-119,ferenczi-holton-zamboni-jam-89,gmp-jtnb-15}),
that a 3iet word is degenerated if and only if the orbits of the
two discontinuity points of the corresponding 3iet $T$ have a
non-empty intersection, formally, $\{T^n(\alpha)\ |\ n\in
\mathbb{Z}\}\cap \{T^n(\alpha+\beta)\ |\ n\in
\mathbb{Z}\}\neq\emptyset$. The complexity of a degenerate 3iet
word is ${\cal C}(n) = n+ const$ for sufficiently large $n$.
Cassaigne~\cite{cassaigne-dlt3} calls one-sided infinite words
with such complexity quasi-sturmian words. By a slight
modification of his results one can show that for any 3iet word
$u$ with complexity $\Ccal_u(n) \leq n+\text{const}$ there exists
a sturmian word $(v_n)_{n\in\Z}$ over $\{0,1\}$ and finite words
$w_1,w_2\in\{A,B,C\}^*$ such that
\[
u = \cdots w_{v_{-2}}w_{v_{-1}}\:|\: w_{v_0} w_{v_1} w_{v_2} \cdots
\,,
\]
that is, $u$ is obtained from $v$ by applying the morphism $0\mapsto w_0$ and
$1\mapsto w_1$.
\end{rem}

\section{Morphisms preserving 3iet words}

\begin{defi}
  A morphism on the alphabet $\{A,B,C\}$ is said to be \emph{3iet preserving} if
  $\varphi(u)$ is a 3iet word for every 3iet word $u$.
\end{defi}

Let us recall that 3iet words are defined as those words coding
3-interval exchange transformations, which are aperiodic.
Similarly, sturmian words are aperiodic words coding 2-interval
exchange transformations.

In the rest of this section we give several useful examples of
3iet preserving morphisms.

\begin{ex}\label{ex:proof-3iet-pres}
  We will prove that the morphism $\varphi$ over $\{A,B,C\}$ given by prescriptions
  \begin{equation}\label{eq:pr-phi}
    A \mapsto AC\,, \qquad
    B \mapsto BC\,, \qquad
    C \mapsto C\,,
  \end{equation}
  is 3iet preserving. Let us consider an arbitrary 3iet word $u$ with
  arbitrary parameters $\alpha,\beta,\gamma$ and $x_0$. The corresponding transformation
  $T$ is given by~\eqref{eq:3iet}. We show that the infinite word
  $\varphi(u)$ is a 3iet word, namely the one with parameters
  $\alpha'=\alpha,\beta'=\beta,\gamma'=\alpha+\beta+\gamma$
  and $x_0'=x_0$.
  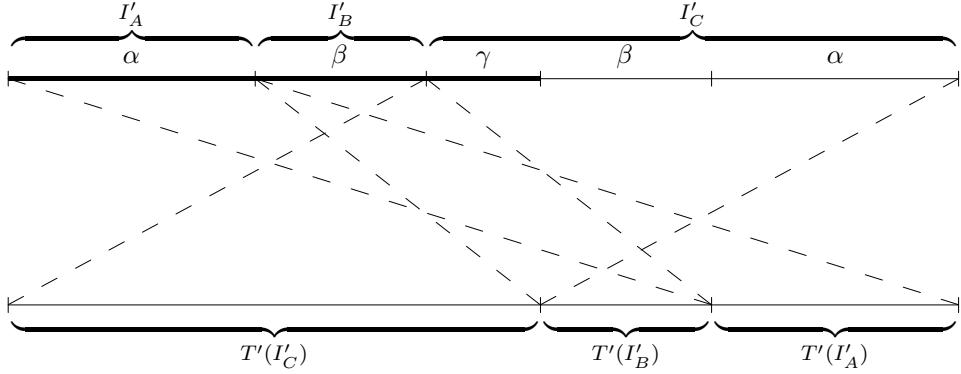
\begin{figure}[ht!]
    \centering
    \setlength{\unitlength}{1cm}
    \begin{picture}(14.5,4.6)(0,0.4)
      \drawline(1,4)(13.5,4)
      \drawline(1,3.9)(1,4.1) \drawline(4.25,3.9)(4.25,4.1) \drawline(6.5,3.9)(6.5,4.1)
      \drawline(8,3.9)(8,4.1) \drawline(10.25,3.9)(10.25,4.1)
      \drawline(13.5,3.9)(13.5,4.1)
      \drawline(1,1)(13.5,1)
      \drawline(1,0.9)(1,1.1) \drawline(8,0.9)(8,1.1)
      \drawline(10.25,0.9)(10.25,1.1) \drawline(13.5,0.9)(13.5,1.1)
      \dashline{0.25}(1,4)(10.25,1) \dashline{0.25}(4.25,4)(13.5,1)
      \dashline{0.25}(4.25,4)(8,1) \dashline{0.25}(6.5,4)(10.25,1)
      \dashline{0.25}(6.5,4)(1,1) \dashline{0.25}(13.5,4)(8,1)
      \put(2.51,4.2){\text{\small$\alpha$}}
      \put(5.25,4.2){\text{\small$\beta$}}
      \put(7.15,4.2){\text{\small$\gamma$}}
      \put(9,4.2){\text{\small$\beta$}}
      \put(11.77,4.2){\text{\small$\alpha$}}
      \put(1.05,4.4){$\overbrace{\hspace{3.15cm}}^{I_A'}$}
      \put(4.3,4.4){$\overbrace{\hspace{2.15cm}}^{I_B'}$}
      \put(6.55,4.4){$\overbrace{\hspace{6.9cm}}^{I_C'}$}
      \put(1.05,0.8){$\underbrace{\hspace{6.9cm}}_{T'(I_C')}$}
      \put(8.05,0.8){$\underbrace{\hspace{2.15cm}}_{T'(I_B')}$}
      \put(10.3,0.8){$\underbrace{\hspace{3.15cm}}_{T'(I_A')}$}
      \linethickness{0.5mm}
      \put(1,4){\line(1,0){7}}
    \end{picture}
    \caption{Graph of the transformation $T'$.}\label{fig:tp-ex}
  \end{figure}

  The transformation $T'$ (see Figure~\ref{fig:tp-ex}) corresponding
  to the parameters $\alpha',\beta',\gamma'$ is given by
  \begin{equation}\label{eq:ex-tp}
    T'(x) =
    \begin{cases}
      x + \alpha+2\beta+\gamma & \text{if $x\in[0,\alpha)\eqcolon I_A'$,} \\
      x +\beta+\gamma & \text{if $x\in[\alpha,\alpha+\beta)\eqcolon I_B'$,} \\
      x -\alpha-\beta & \text{if $x\in[\alpha+\beta,2\alpha+2\beta+\gamma)\eqcolon I_C'$.}
    \end{cases}
  \end{equation}

  Obviously, \eqref{eq:3iet} and~\eqref{eq:ex-tp}
  imply for a point $x\in I_A=I_A'$ that
  \begin{align*}
    T'(x) & = x + \alpha + 2\beta + \gamma \in I_C'\,, \\
    (T')^2(x) & = x + \beta + \gamma = T(x)\,.
  \end{align*}
  Hence any point $x\in I_A$ belongs in the new 3iet to the interval $I_A'$,
  its first iteration is $T'(x)\in I_C'$ and the second iteration
  $(T')^2(x)$ sends to the same place as the first iteration of the
  original transformation $T$.  Therefore we substitute $A\mapsto AC$.
  Similarly, for a point $x\in I_B=I_B'$ we have
  \begin{align*}
    T'(x) & = x + \beta + \gamma \in I_C'\,, \\
    (T')^2(x) & = x - \alpha + \gamma = T(x)\,,
  \end{align*}
  and so  $B\mapsto BC$.
  Finally, for $x\in I_C \subsetneqq I_C'$ we get $T'(x) = T(x)$ and hence $C\mapsto C$.
 Thus we see that the 3iet word coding $x'_0$ under $T'$ coincides with
 the word $\varphi(u)$.
\end{ex}

\begin{ex}\label{ex:prehod-AC}
  It is easy to see that the morphism $\xi$ over $\{A,B,C\}$ given by prescriptions
  \begin{equation}\label{eq:pr-xi}
    A \mapsto C\,, \qquad B \mapsto B\,, \qquad C \mapsto A\,,
  \end{equation}
  is a 3iet preserving morphism. To a 3iet word, which codes the orbit of $x_0$
  under the transformation $T$ with intervals
  $[0,\alpha)\cup[\alpha,\alpha+\beta)\cup[\alpha+\beta,\alpha+\beta+\gamma)$,
  is assigns a 3iet word, which codes the orbit of $\alpha+\beta+\gamma-x_0$
  under the transformation $\tilde{T}$ with intervals
  $(0,\gamma]\cup(\gamma,\gamma+\beta]\cup(\gamma+\beta,\gamma+\beta+\alpha]$.
\end{ex}

\begin{ex}\label{ex:primitive-3iet-pres}
  Let us consider the morphism $\varphi_0$ on $\{A,B,C\}$ given by $A\mapsto B$,
  $B\mapsto BCB$ and $C\mapsto CAC$. It is a primitive morphism with
  $\det\mat{M}_{\varphi_0}=1$ and $\mat{M}_{\varphi_0}^3>0$. Let $u$ be an arbitrary
  3iet word with parameters $\alpha,\beta,\gamma$. Using the same technique as in
  Example~\ref{ex:proof-3iet-pres} one can show that $\varphi_0$ is 3iet preserving;
  the 3iet word coinciding with $\varphi_0(u)$ has parameters $\alpha'=\gamma$,
  $\beta'=\beta+\alpha+\beta$, $\gamma'=\gamma+\beta+\gamma$.
\end{ex}

\section{Proof of Theorem A}

The aim of this section is to prove that the matrix $\mat{M}$
of a 3iet preserving morphism fulfills the following condition
\begin{equation}\label{eq:meme-cond}
  \mat{M}\mat{E}\mat{M}^T = \pm\mat{E}, \qquad\text{where }
  \mat{E} =
  \Bigl(\!\begin{smallmatrix}
    0 & 1 & 1 \\
    -1 & 0 & 1 \\
    -1 & -1 & 0
  \end{smallmatrix}\Bigr)\,.
\end{equation}

The main tool used in the proof of this property of $\mat{M}$ is
the fact that the matrix of a sturmian morphism has determinant
$\pm1$ and some auxiliary statements formulated as
Lemma~\ref{lem:podprostorP} and Lemma~\ref{lem:e1e2e3}.

\begin{lem}\label{lem:podprostorP}
Let $\varphi$ be a 3iet preserving morphism and $\mat{M}$ its
incidence matrix. Let ${\cal P}$ be a subspace of $\R^3$ spanned
by the vectors
$\Bigl(\begin{smallmatrix}1\\1\\0\end{smallmatrix}\Bigr)$,
$\Bigl(\begin{smallmatrix}0\\1\\1\end{smallmatrix}\Bigr)$. Then
$\mat{M}{\cal P}={\cal P}$.
\end{lem}

\begin{proof}
If $u$ is a 3iet word with parameters $(\alpha,\beta,\gamma)$, then
according to~(\ref{eq:density-phi-u}) and
Corollary~\ref{coro:hustoty}, $\varphi(u)$ is a 3iet word with
parameters $(\alpha,\beta,\gamma)\mat{M}$. Since $\varphi$ is a 3iet
preserving morphism, it means that $\varphi(u)$ is aperiodic,
whenever $u$ is aperiodic. With the help of
Corollary~\ref{coro:period}, it implies that for every pair
$K,L\in\Z\setminus\{0\}$ and every triple of positive numbers
$(\alpha,\beta,\gamma)$, we have
\begin{equation}\label{eq:zachovava}
(\alpha,\beta,\gamma)\mat{M}\!\!\left(\!\!\begin{smallmatrix}K\\K+L\\L\end{smallmatrix}\!\!\right)
= 0 \quad\implies\quad \exists\, H,S\in\Z\setminus\{0\} \ \hbox{ such that }\
(\alpha,\beta,\gamma)\left(\!\!\begin{smallmatrix}H\\H+S\\S\end{smallmatrix}\!\!\right)
= 0\,.
\end{equation}
Since
$\Bigl\{\mat{M}\!\!\left(\!\!\begin{smallmatrix}K\\K+L\\L\end{smallmatrix}\!\!\right)
\,\Bigm|\, K,L\in\Z \Bigr\}$ is a 2-dimensional lattice in $\R^3$,
there exist two linearly independent pairs $K_1,L_1$, $K_2,L_2$
such that
$\mat{M}\!\!\left(\!\!\begin{smallmatrix}K_i\\K_i+L_i\\L_i\end{smallmatrix}\!\!\right)$,
$i=1,2$, have both positive and negative components, and therefore
for both $i=1,2$, there exist infinitely many triples
$(\alpha,\beta,\gamma)$ such that
$(\alpha,\beta,\gamma)\mat{M}\!\!\left(\!\!\begin{smallmatrix}
K_i\\K_i+L_i\\L_i\end{smallmatrix}\!\!\right)=0$. This, together
with~\eqref{eq:zachovava}, implies
\begin{equation}\label{eq:901}
\mat{M}\!\!\left(\!\!\begin{smallmatrix}
K_i\\K_i+L_i\\L_i\end{smallmatrix}\!\!\right) = {\it
const.}\!\!\left(\!\!\begin{smallmatrix}
H_i\\H_i+S_i\\S_i\end{smallmatrix}\!\!\right)\,,\quad \hbox{ for
some } \ H_i,S_i\in\Z\setminus\{0\}, \ i=1,2.
\end{equation}
Consequently, $\mat{M}{\cal P}\subseteq{\cal P}$. We now show that
$\mat{M}{\cal P}={\cal P}$. Suppose the opposite, i.e., that
$\mat{M}\!\!\left(\begin{smallmatrix}
1\\1\\0\end{smallmatrix}\right)$ and
$\mat{M}\!\!\left(\begin{smallmatrix}
0\\1\\1\end{smallmatrix}\right)$ are linearly dependent. Then
there exist $K,L\in\Z\setminus\{0\}$ such that
$$
\left(\begin{smallmatrix} 0\\0\\0\end{smallmatrix}\right)=
K\mat{M}\left(\begin{smallmatrix} 1\\1\\0\end{smallmatrix}\right)+
L \mat{M}\left(\begin{smallmatrix}
0\\1\\1\end{smallmatrix}\right)=\mat{M}\left(\!\!\begin{smallmatrix}
K\\K+L\\L\end{smallmatrix}\!\!\right)\,.
$$
This, however, implies that for arbitrary parameters
$(\alpha,\beta,\gamma)$, we have
$$
(\alpha,\beta,\gamma)\mat{M}\left(\!\!\begin{smallmatrix}
K\\K+L\\L\end{smallmatrix}\!\!\right) = 0\,,
$$
i.e., the word $\varphi(u)$ is periodic for arbitrary 3iet word $u$, which
is a contradiction with the assumption that $\varphi$ is a 3iet
preserving morphism.
\end{proof}

\begin{rem}\label{ex:regular}
Denote
$$
\vec{x}_1\coloneq
\Bigl(\begin{smallmatrix}1\\1\\0\end{smallmatrix}\Bigr)\,,\qquad
\vec{x}_2\coloneq
\Bigl(\begin{smallmatrix}0\\1\\1\end{smallmatrix}\Bigr)\,,\qquad
\vec{x}_3\coloneq
\Bigl(\begin{smallmatrix}0\\1\\0\end{smallmatrix}\Bigr)\,.
$$
The triplet of vectors $\vec{x}_1$, $\vec{x}_2$, $\vec{x}_3$ forms
a basis of $\R^3$. Denoting
$\mat{P}=\Bigl(\begin{smallmatrix}1&0&0\\1&1&1\\0&1&0\end{smallmatrix}\Bigr)\,,$
we have $\det\mat{P}=1$, and thus $\vec{x}_1$, $\vec{x}_2$,
$\vec{x}_3$ is also a basis of the integer lattice $\Z^3$. In the
same time, the pair $\vec{x}_1$, $\vec{x}_2$ is a basis of the
invariant subspace ${\cal P}$ of the matrix $\mat{M}$. We have
\[
\mat{P}^{-1}=\Bigl(\!\begin{smallmatrix}1&0&0\\0&0&1\\-1&1&-1\end{smallmatrix}\Bigr)\,
\quad\text{ and }\quad
\mat{P}^{-1}\mat{M}\mat{P} =
\begin{pmatrix}
  m_{11}+m_{12} & m_{12}+m_{13} & m_{12}\\
  m_{31}+m_{32} & m_{32}+m_{33} & m_{32}\\
  0&0& -m_{12} +m_{22} - m_{32}
\end{pmatrix}\,,
\]
 where the 0's in the third row correspond to the fact that
$\mat{P}^{-1}\mat{M}\mat{P}$ can be seen as the matrix $\mat{M}$
written in the basis $\vec{x}_1$, $\vec{x}_2$, $\vec{x}_3$, where
the first two vectors form a basis of the invariant subspace
${\cal P}$. Since $\mat{M}{\cal P}={\cal P}$, we have
\[
  \det\begin{pmatrix}
  m_{11}+m_{12} & m_{12}+m_{13}\\
  m_{31}+m_{32} & m_{32}+m_{33}
  \end{pmatrix} \neq 0\,.
  \]
\end{rem}

\begin{lem}\label{lem:e1e2e3}
  Let $\mat{M}=(m_{ij})$ be the incidence matrix of a 3iet preserving morphism $\varphi$.
  Then
  \begin{equation}\label{eq:mala-delta}
  \det\begin{pmatrix}
  m_{11}+m_{12} & m_{12}+m_{13}\\
  m_{31}+m_{32} & m_{32}+m_{33}
  \end{pmatrix} = \delta \in \{1,-1\}\,.
  \end{equation}
\end{lem}

\begin{proof}
  Let us choose a sturmian word $u\in\{A,C\}^\Z$ and a sequence $(u^{(m)})_{m\in\N}$ of 3iet
  words  such that $u=\lim_{m\rightarrow\infty} u^{(m)}$.
  For example, let $u$ be the coding of
  $x_0=0$ under the 2-interval exchange transformation $T$ with $I_0=[0,1-\alpha)$
  and $I_1=[1-\alpha,1)$, where $\alpha$ is an arbitrary irrational number. Then we can choose
  $u^{(m)}$ to be the 3iet word that codes $x_0=0$ under the 3-interval exchange
  transformation with intervals $I_A=[0,1-\alpha - \frac{1}{m})$,
  $I_B=[1-\alpha -\frac{1}{m},1-\alpha)$ and $I_C=[1-\alpha,1)$.

  Let $\sigma$ be a morphism given by
  \[
  A \mapsto A\,,\quad B \mapsto AC\,, \quad C \mapsto C\,.
  \]
  Since any morphism on $\{A,B,C\}^\Z$ is a continuous mapping, we have
  \[
  (\sigma\circ\varphi)(u^{(m)}) \rightarrow (\sigma\circ\varphi)(u)\,.
  \]
  According to the assumption, the morphism $\varphi$ is 3iet preserving, hence
  $\varphi(u^{(m)})$ are 3iet words. By
  Lemma~\ref{lem:delta-of-3iet},
  the words $(\sigma\circ\varphi)(u^{(m)})$, $m\in\N$, code 2-interval exchange transformations,
  and by Lemma~\ref{lem:limit-of-sturm}, the limit of these words, that is
  the word $(\sigma\circ\varphi)(u)$, is either sturmian or the densities
  of its letters are rational.

  The matrix of $\sigma$ is
  $\Bigl(\!\begin{smallmatrix} 1 & 0 & 0 \\ 1 & 0 & 1 \\ 0 & 0 & 1\end{smallmatrix}\Bigr)$,
  which implies by~(\ref{eq:matrix-compose}) that the matrix of $\sigma\circ\varphi$
  is
  \[
  \mat{M}_{\sigma\circ\varphi} =
  \begin{pmatrix}
   m_{11}+m_{12} & 0 & m_{12}+m_{13} \\
   m_{21}+m_{22} &  0 & m_{22}+m_{23} \\
   m_{31}+m_{32} &   0 &  m_{32}+m_{33}
  \end{pmatrix} \,.
  \]
  Since $\sigma\circ\varphi$ maps a sturmian word $u$ over $\{A,C\}$ to a word
  over the same alphabet, we are interested only in the matrix of this morphism over
  $\{A,C\}$, that is,
  \begin{equation}\label{eq:msvlnkou}
  \widetilde{\mat{M}}\ = \
  \begin{pmatrix}
    m_{11}+m_{12} & m_{12}+m_{13} \\
    m_{31}+m_{32} & m_{32}+m_{33}
  \end{pmatrix} \,.
  \end{equation}

  Let us suppose that the densities of $A$ and $C$ in $u$ are $1-\alpha$ and $\alpha$,
  respectively. Using~(\ref{eq:density-phi-u}) we find the density of $A$ in
  $(\sigma\circ\varphi)(u)$ to be

\begin{equation}\label{eq:density-b}
  \rho(A) = \frac{(1-\alpha, \alpha)\ \widetilde{\mat{M}}
  \Bigl(\begin{smallmatrix}1\\0\end{smallmatrix}\Bigr)}
  {(1-\alpha, \alpha)\ \widetilde{\mat{M}}
  \Bigl(\begin{smallmatrix}1\\1\end{smallmatrix}\Bigr)}\,.
\end{equation}
  If $\rho(A)$ is irrational, the word $(\sigma\circ\varphi)(u)$ is sturmian
  and hence the morphism $\sigma\circ\varphi$ is sturmian. This implies
  $\det\widetilde{\mat{M}}=\pm 1$.

  The irrational number $\alpha$, i.e., the density of $A$ in the sturmian word $u$,
  was chosen arbitrarily. Therefore $\rho(A)$, given by~(\ref{eq:density-b}),
  will be rational for any irrational $\alpha$ only in case when

 \begin{equation}\label{eq:rational-rhoB}
  p\widetilde{\mat{M}}\Bigl(\begin{smallmatrix}1\\0\end{smallmatrix}\Bigr) = q
  \widetilde{\mat{M}}\Bigl(\begin{smallmatrix}1\\1\end{smallmatrix}\Bigr),
  \quad \hbox{for some}  \ \ p,q \in \Z\setminus\{0\}\,.
  \end{equation}
  This however implies that the matrix $\widetilde{\mat{M}}$ is
  singular, which contradicts Remark~\ref{ex:regular}.
\end{proof}

We are now in position to finish the proof of Theorem A.

\begin{thma}
  Let $\mat{M}$ be the incidence matrix of a 3iet preserving morphism.
  Then
  \begin{equation}\label{eq:meme}
    \mat{M}\mat{E}\mat{M}^T = \pm\mat{E},\quad\text{where }\ \mat{E} =
    \Bigl(\!\begin{smallmatrix}0&1&1\\-1&0&1\\-1&-1&0\end{smallmatrix}\Bigr)\,.
  \end{equation}
\end{thma}
\begin{proof}
Using the notation of Remark~\ref{ex:regular} for the matrix
$\mat{P}$, we obviously see that the matrix
$\mat{P}^{-1}\mat{M}\mat{P}$ has $(0,0,-1)$ for its left
eigenvector corresponding to the eigenvalue $-m_{12} +m_{22} -
m_{32}$. It is then trivial to verify that
$(0,0,-1)\mat{P}^{-1}=(1,-1,1)$ is a left eigenvector of the
matrix $\mat{M}$ corresponding to the same eigenvalue. Since
\begin{equation}\label{eq:det_PMP}
\det{\mat{M}} = \det (\mat{P}^{-1}\mat{M}\mat{P}) = \delta (-m_{12} +m_{22} - m_{32})\,,
\end{equation}
where $\delta\in\{-1,1\}$ is given by~(\ref{eq:mala-delta}),
we derive that $(1,-1,1)$ is a left eigenvector of the matrix $\mat{M}$ corresponding to the
eigenvalue $\delta\det{\mat{M}}$. Denoting $\Delta\coloneq \det{\mat{M}}$, we can write
\begin{equation}\label{eq:vector}
(1,-1,1)\mat{M}= \delta\Delta (1,-1,1)\,.
\end{equation}
This implies that the matrix $\mat{M}$ can be written in the following form,
\begin{equation}\label{eq:matM}
\mat{M}=
\begin{pmatrix}
  m_{11} & m_{12} & m_{13} \\
  m_{11}+m_{31}-\delta\Delta & m_{12}+m_{32}+\delta\Delta & m_{13}+m_{33} -\delta\Delta \\
  m_{31}& m_{32} & m_{33}
\end{pmatrix}\,.
\end{equation}
With this, one can verify by inspection, that $\mat{M}\mat{E}\mat{M}^T = \delta\mat{E}$,
using Lemma~\ref{lem:e1e2e3} for simplification of algebraic expressions.
\end{proof}

As a partial result, we have shown in the above proof the following interesting statement.

\begin{coro}\label{coro:eigenvalues-M}
  Let $\mat{M}$ be the matrix of a 3iet preserving morphism $\varphi$.
  Then the vector $(1,-1,1)$ is a left eigenvector
  of $\mat{M}$, associated with the eigenvalue $\det\mat{M}$ or $-\det\mat{M}$, i.e.,
  \begin{equation}\label{eq:left-eigen}
    (1,-1,1)\mat{M} = \pm\det\mat{M}(1,-1,1)\,.
  \end{equation}
The other eigenvalues $\lambda_1$ and $\lambda_2$ of the matrix
$\mat{M}$ are either quadratic mutually conjugate algebraic units,
or $\lambda_1, \lambda_2 \in \{1, -1\}$.
\end{coro}

From the form~(\ref{eq:matM}) of the matrix $\mat{M}$ we derive
the following Corollary.

\begin{coro}\label{coro:sumy-radku}
Let $\mat{M}$ be a matrix of a 3iet preserving morphism. Then the
sum of its first and the third row differs from the sum of its
second row by $\pm\det\mat{M}$. Formally,
\[
(1,0,1)\mat{M}\Big(\begin{smallmatrix}1\\1\\1\end{smallmatrix}\Big) -
(0,1,0)\mat{M}\Big(\begin{smallmatrix}1\\1\\1\end{smallmatrix}\Big) = \pm\det\mat{M}\,.
\]
\end{coro}

\section{3iet preserving morphisms versus fixed points}

The proof of Theorem B, which is performed in
Section~\ref{sec:proof-B} is based on the properties of 3iet
words, which are fixed points of morphisms. In this section we
therefore inspect, which 3iet preserving morphisms have a fixed
point.

A fixed point of a morphism $\varphi$ over an alphabet $\Acal$ is
the limit $\lim_{n\rightarrow\infty}\varphi^n(a_i)|\varphi^n(a_j)$
for some $a_i,a_j\in\Acal$. Similarly to the case of sturmian
words, the set of 3iet words is not compact, and therefore in
general the accumulation point $u$ of a sequence
$(u^{(m)})_{m\in\N}$ of 3iet words is not necessarily a 3iet word.
The special case when the accumulation point belongs to the set of
3iet words is treated by the following Lemma.

\begin{lem}\label{lem:apendix}
  Let $\alpha,\beta,\gamma$ be positive real numbers such that
  $\alpha+\beta$ and $\beta+\gamma$ are linearly independent over $\Q$.
  Let $T_1$, $T_2$ be the 3iet transformations with parameters
  $\alpha,\beta,\gamma$ and domain $[0,\alpha+\beta+\gamma)$,
    $(0,\alpha+\beta+\gamma]$, respectively. Let $(u^{(n)})_{n\in\N}$
  be a sequence of 3iet words and $(x^{(n)})_{n\in\N}$ a sequence of
  points in $[0,\alpha+\beta+\gamma]$ such that
  \begin{itemize}
  \item
    $u^{(n)}=u_{T_1}(x^{(n)})$ or $u^{(n)}=u_{T_2}(x^{(n)})$ for all $n\in\N$;
  \item
    $x^{(n)}$ is a monotonous sequence with the limit $x$.
  \end{itemize}
  Then $\lim_{n\to\infty}u^{(n)}$ exists and is equal to the 3iet
  word $u_{T_1}(x)$ or $u_{T_2}(x)$.
\end{lem}
\begin{proof}
We use a statement from~\cite{gmp-jtnb-15}.
For a given $m$ put
$$
D_m\coloneq \{T^{i}_1(\alpha), T^{i}_1(\alpha+\beta),
T^{i}_2(\alpha), T^{i}_2(\alpha+\beta)\mid -m\leq i\leq m \}\,.
$$
Let $a<b$ and let $(a,b)\cap D_m=\emptyset$. Then for all
$z\in(a,b)$ we have
\begin{eqnarray}
\label{eq:Dmi}
\dd\bigl(u_{T_1}(a),u_{T_1}(z)\bigr) < \frac{1}{1+m}\,, &\quad&
\dd\bigl(u_{T_1}(a),u_{T_2}(z)\bigr) < \frac{1}{1+m}\,,\\[2mm]
\dd\bigl(u_{T_2}(b),u_{T_1}(z)\bigr) < \frac{1}{1+m}\,, &\quad&
\dd\bigl(u_{T_2}(b),u_{T_2}(z)\bigr) < \frac{1}{1+m}\,.
\label{eq:Dmii}
\end{eqnarray}
Assume that the sequence $(x^{(n)})_{n\in\N}$ is decreasing. For
$\varepsilon>0$, we find $m\in\N$ such that
$\varepsilon>\frac1{m+1}$ and we put $\delta_m\eqcolon \sup\{y>x
\mid y\notin D_m \}$. Since $x^{(n)}\searrow x$, there exists
$n_0$ such that for all $n>n_0$ we have $x\leq
x^{(n)}<x+\delta_m$. Since $u^{(n)}=u_{T_1}(x^{(n)})$ or
$u^{(n)}=u_{T_2}(x^{(n)})$, we obtain, using~\eqref{eq:Dmi} for
the interval $(a,b)=(x,x+\delta_m)$, that
$\dd\bigl(u_{T_1}(x),u^{(n)}\bigr)<\varepsilon$, which implies
$\lim_{n\to\infty} u^{(n)}=u_{T_1}(x)$. Similarly we
use~\eqref{eq:Dmii} in the case $x^{(n)}\nearrow x$.
\end{proof}

\begin{rem}
  The assumption of primitivity of the morphism $\varphi$ is essential in the above statement.
  For example, the morphism $\varphi$ defined by~\eqref{eq:pr-phi} is 3iet preserving,
  yet the only fixed points of an arbitrary power $\varphi^p$, $p\in\N$, $p\geq 1$, are
  \[
  \cdots CCC|ACCC \cdots, \qquad
  \cdots CCC|BCCC \cdots, \qquad
  \cdots CCC|CCC \cdots
  \]
\end{rem}

The following proposition deals with the original aim of this
section, namely with the search for 3iet preserving morphisms
having 3iet words as their fixed points.

\begin{prop}\label{prop:fixedpoint}
  Let $\varphi$ be a primitive 3iet preserving morphism.
  Then there exists $p\in\N$, $p\geq 1$, such
  that $\varphi^p$ has a fixed point, and this fixed point is a 3iet word.
\end{prop}
\begin{proof}
Without loss of generality, we may assume that the incidence
matrix $\mat{M}$ of the morphism $\varphi$ is positive. Otherwise,
we show the validity of the statement for $\psi=\varphi^k$ for
some $k$, which implies the validity of the statement for
$\varphi$.

Let $(\alpha,\beta,\gamma)$ be a positive left eigenvector of
$\mat{M}$. First we show that an infinite word coding a 3iet with
such parameters is not periodic. For contradiction, assume that
$(\alpha,\beta,\gamma)$ satisfy~\eqref{coro:period}, that is,
\begin{equation}\label{eq:xxxx}
(\alpha,\beta,\gamma)\left(\!\begin{smallmatrix}K\\K+L\\L\end{smallmatrix}\!\right)
= 0\,, \quad\text{for some $K,L\in\Z\setminus\{0\}$.}
\end{equation}
If the Perron eigenvalue $\lambda_1$ of $\mat{M}$ is a quadratic
irrational number, one can assume without loss of generality that
the components of the vector $(\alpha,\beta,\gamma)$ belong to the
quadratic field $\mathbb{Q}(\lambda_1)$. For any $x\in
\mathbb{Q}(\lambda_1)$, denote by $x'$ the image of $x$ under the
Galois automorphism of $\mathbb{Q}(\lambda_1)$. Since the matrix
$\mat{M}$ and the vector
$\left(\!\begin{smallmatrix}K\\K+L\\L\end{smallmatrix}\!\right)$
have integer components, the vector $(\alpha',\beta',\gamma')$ is
an eigenvector to the eigenvalue $\lambda'_1 = \lambda_2$ and
satisfies
\[
(\alpha',\beta',\gamma')\left(\!\begin{smallmatrix}K\\K+L\\L\end{smallmatrix}\!\right)
= 0\,.
\]
Using Corollary~\ref{coro:eigenvalues-M}, the vector $(1, -1,1)$
is a left eigenvector of $\mat{M}$ corresponding to the eigenvalue
$\pm \det\mat{M}$. Therefore vectors $(\alpha',\beta',\gamma')$,
$(\alpha,\beta,\gamma)$ and $(1, -1,1)$ are eigenvectors of
$\mat{M}$ corresponding to different eigenvalues, which means that
they are linearly independent. All of them are orthogonal to the
vector
$\left(\!\begin{smallmatrix}K\\K+L\\L\end{smallmatrix}\!\right)$,
which implies $K=L=0$. This contradicts~\eqref{eq:xxxx}.

By Corollary~\ref{coro:eigenvalues-M}, it remains to discuss the
case when the Perron eigenvalue of $\mat{M}$ is $\lambda_1=1$.
This is impossible due to the fact that a positive integral matrix
$\mat{M}$ cannot have $1$ as its eigenvalue corresponding to a
positive eigenvector. Thus we have shown that the infinite word
coding a 3iet with parameters $\alpha,\beta,\gamma$ is not
periodic.

Denote $T_1$, $T_2$ the 3iet transformations with parameters
$\alpha,\beta,\gamma$ and domain $[0,\alpha+\beta+\gamma)$,
$(0,\alpha+\beta+\gamma]$, respectively.

Let $u^{(0)}$ be an arbitrary 3iet word coding the orbit of a
point by $T_1$. Put
\[
u^{(n)}\coloneq \varphi^n(u^{(0)})\,, \quad\text{for $n\geq 1$.}
\]
Since the vector of densities of $u^{(0)}$ is a left eigenvector
of the incidence matrix of the morphism $\varphi$, every word
$u^{(n)}$, $n\in\N$, has the same density of letters. As $\varphi$
is a 3iet preserving morphism, the word $u^{(n)}$ is a 3iet word
coding the orbit of a point under $T_1$ or $T_2$, for every
$n\in\N$.

The space of infinite words over the alphabet $\{A,B,C\}$ is
compact, and thus there exists a Cauchy subsequence of the
sequence $(u^{(n)})_{n\in\N}$. Therefore there exist
$m_0,n_0\in\N$, $n_0>m_0$, such that
\begin{equation}\label{eq:401}
\dd\bigl(u^{(n_0)},u^{(m_0)}\bigr)<\frac12\,.
\end{equation}
Set $p\coloneq n_0-m_0$ and $v=\cdots v_{-2}v_{-1}|v_0v_1\cdots
\coloneq u^{(m_0)}$. Since $ u^{(n_0)} =
\varphi^{n_0-m_0}(u^{(m_0)}) = \varphi^p(v)$,
inequality~\eqref{eq:401} can be rewritten as
\begin{equation}\label{eq:402}
\dd\bigl(\varphi^p(v),v\bigr)<\frac12 \,.
\end{equation}
The latter, together with the primitivity of the morphism
$\varphi$, implies
\[
\varphi^p(v_0)=v_0w_0\quad\hbox{ and }\quad
\varphi^p(v_{-1})=w_{-1}v_{-1}
\]
for some non-empty words $w_0,w_{-1}\in\{A,B,C\}^*$. Therefore the
morphism $\varphi^p$ has the fixed point
\[
\lim_{n\to\infty} \varphi^{np}(v)\,.
\]
Since $\varphi^{np}(v)$ is a 3iet word given by $T_1$, or $T_2$,
there exist for every $n$ a number
$x^{(n)}\in[0,\alpha+\beta+\gamma]$, such that
\[
\varphi^{np}(v) = u_{T_1}(x^{(n)})\ \hbox{ or }\
u_{T_2}(x^{(n)})\,.
\]
Denote by $x$ the limit of some monotonous subsequence of
$(x^{(n)})_{n\in\N}$, i.e., $ x=\lim_{n\to\infty}x^{(k_n)}$.
According to Lemma~\ref{lem:apendix},
\[
\lim_{n\to\infty} \varphi^{np}(v) = u_{T_1}(x)\ \hbox{ or }\
u_{T_2}(x)\,,
\]
which means that $\varphi^p$ has as its fixed point a 3iet word,
namely $u_{T_1}(x)$ or $u_{T_2}(x)$, respectively.
\end{proof}

\section{Proof of Theorem~B}\label{sec:proof-B}

In the proof of Theorem B we use certain properties of discrete
sets associated with 3iet words. Every 3iet word can be
geometrically represented using a C\&P set. On the other hand,
every fixed point of a primitive morphism can be represented by a
self-similar set, which is constructed using a right eigenvector
of the matrix of the morphism. The crucial point in the proof of
Theorem B is the fact that for a 3iet word being a fixed point of
a primitive morphism these two geometric representations coincide.

We first show that the determinant of the incidence matrix of a
3iet preserving morphism is in modulus smaller than 1. For that we
use the following technical lemma.

\begin{lem}\label{lem:Pk}
Let $\varepsilon\in(0,1)$ be a quadratic irrational number with
conjugate $\varepsilon'<0$. Let $\lambda\in(0,1)$ be a quadratic
unit such that its conjugate satisfies $\lambda'>1$ and
$\lambda'\Z[\varepsilon']=\Z[\varepsilon'] \coloneq \Z +
\varepsilon'\Z$. Let us denote $\Lambda\coloneq\lambda'$,
$\eta\coloneq-\varepsilon'$ and
\[
P_n(x)\coloneq
\#\bigl(x,x+(1+2\eta)\Lambda^n\bigr]\cap\Sigma_{\varepsilon,\eta}(\Omega)\,,
\]
where $\Omega$ is a bounded interval. Then there is a constant $R$ such that
\[
|P_n(x) - P_n(y)| \leq R\,,
\]
for any $x,y\in\R$ and $n\in\N$.
\end{lem}
The proof of Lemma exploits some simple properties of C\&P sets,
which are however not related to infinite words. Therefore we
postpone it to the appendix.

\begin{prop}\label{prop:detM-01}
The incidence matrix $\mat{M}$ of a primitive 3iet preserving
morphism $\varphi$ satisfies $|\det\mat{M}|\leq 1$.
\end{prop}

\begin{proof}
Without loss of generality we assume that $\varphi$ has a 3iet
fixed point $u$, and, moreover, that both the matrix $\mat{M}$ and
its spectrum are positive. This is possible since according to
Proposition~\ref{prop:fixedpoint} for any primitive 3iet
preserving morphism there exists $p\in\N$ such that $\varphi^p$
has a fixed point and $|\det\mat{M}|\leq 1\ \Leftrightarrow\
|\det\mat{M}^p|\leq 1$.

Let us denote by $\Lambda$ the dominant (Perron) eigenvalue of
$\mat{M}$. Its second eigenvalue is by
Corollary~\ref{coro:eigenvalues-M} equal to $\pm\det\mat{M}$, the
third one is denoted by $\lambda$. A positive integer matrix
cannot have $1$ as its dominant eigenvalue, hence by
Corollary~\ref{coro:eigenvalues-M}, $\Lambda>1$ is a quadratic
algebraic unit such that $\Lambda'=\lambda$.

Without loss of generality we assume that a positive right
eigenvector associated with the Perron eigenvalue $\Lambda$ is
such that the modulus of its third component is greater than the
modulus of the first one. Otherwise, we use
$\xi\circ\varphi\circ\xi$ instead of $\varphi$, where $\xi$ is
defined as in Example~\ref{ex:prehod-AC}. Matrices corresponding
to $\varphi$ and $\xi\circ\varphi\circ\xi$ have the same spectrum,
the first and the last component of eigenvectors being
interchanged.

The fixed point $u$ of $\varphi$ is the coding of a 3iet with parameters
$\alpha,\beta,\gamma$, with a starting point $x_0$. By~(\ref{eq:density-phi-u})
and Corollary~\ref{coro:hustoty} the vector $(\alpha,\beta,\gamma)$ is a left
eigenvector of $\mat{M}$ corresponding to $\Lambda$.

In the proof we use properties of a C\&P set; we construct it in
such a way that it coincides with the geometric representation of
the fixed point $u$. Let us define parameters $\varepsilon,l,c$
and the interval $\Omega$ by~(\ref{eq:prevod}). Note that
$(l-\varepsilon,1-l,l-1+\varepsilon)$ is also an eigenvector to
$\Lambda$, because it is just a scalar multiple of
$(\alpha,\beta,\gamma)$, and, moreover, since $\Lambda$ is a
quadratic irrational number, the parameters $\varepsilon, l$
belong to the same quadratic algebraic field $\Q(\Lambda)$. By
$x'$ we denote the image of $x\in\Q(\Lambda)$ under the Galois
automorphism on $\Q(\Lambda)$.

Let us denote $\vec{F} =
\Big(\!\begin{smallmatrix}-\varepsilon\\1-2\varepsilon\\1-\varepsilon\end{smallmatrix}\!\Big)$.
The vector $\vec{F}$ is orthogonal to two left eigenvectors $(1,-1,1)$ and
$(l-\varepsilon,1-l,l-1+\varepsilon)$ associated with eigenvalues
$\pm\det\mat{M}$ and $\Lambda$, respectively. The matrix $\mat{M}$ has three
different eigenvalues, therefore $\vec{F}$ is a right eigenvector to the third
eigenvalue $\lambda$.

Since the matrix $\mat{M}$ is integral, the vector $\vec{F}'
\coloneq
\bigg(\!\begin{smallmatrix}-\varepsilon'\\1-2\varepsilon'\\1-\varepsilon'\end{smallmatrix}\!\bigg)$
is a right eigenvector corresponding to the dominant eigenvalue
$\lambda'=\Lambda$, that is,
\begin{equation}\label{eq:MF-ML}
\mat{M}
\bigg(\!\begin{smallmatrix}-\varepsilon'\\1-2\varepsilon'\\1-\varepsilon'\end{smallmatrix}\!\bigg)
= \lambda'
\bigg(\!\begin{smallmatrix}-\varepsilon'\\1-2\varepsilon'\\1-\varepsilon'\end{smallmatrix}\!\bigg)
= \Lambda
\bigg(\!\begin{smallmatrix}-\varepsilon'\\1-2\varepsilon'\\1-\varepsilon'\end{smallmatrix}\!\bigg)\,.
\end{equation}
Therefore the components of $\vec{F}'$ are either all positive or
all negative. By assumption, the modulus of the third component of
a right dominant eigenvector is greater than the modulus of the
first one, which implies that all components of $\vec{F}'$ are
positive, i.e., $-\varepsilon'>0$.

We define a C\&P set with parameters $\varepsilon,\eta,\Omega$,
where $\varepsilon$ and $\Omega$ are as above and we put
$\eta\coloneq-\varepsilon'$. By Proposition~\ref{prop:cap}, the
distances between adjacent elements of
$\Sigma_{\varepsilon,\eta}(\Omega)$ take values $\mu_A=\eta$,
$\mu_B=1+2\eta$, and $\mu_C = 1 + \eta$ and their ordering with
respect to the origin is coded by the word $u$. Let
$(t_n)_{n\in\Z}$ denote a strictly increasing sequence such that
$\Sigma_{\varepsilon,\eta}(\Omega) = \{t_n\:|\:n\in\Z\}$.
According to Section~\ref{sec:geometric}, this C\&P set is also
the geometric representation of the fixed point $u$ of $\varphi$
and
\[
\#\big((\Lambda t_n,\Lambda t_{n+1}]\cap \Sigma_{\varepsilon,\eta}(\Omega)\big) =
\begin{cases}
|\varphi(A)| & \text{if $t_{n+1}-t_n = \mu_A = \eta\,,$} \\
|\varphi(B)| & \text{if $t_{n+1}-t_n = \mu_B = 1+2\eta\,.$} \\
|\varphi(B)| & \text{if $t_{n+1}-t_n = \mu_C = 1+\eta\,.$}
\end{cases}
\]

As the fixed point of a morphism is also the fixed point of an arbitrary power
of this morphism, the geometric representations of $\varphi$ and $\varphi^n$ are the same
for any $n\in\N$. Since $AC$ is a factor of any 3iet word, there exist $k,m\in\N$ such that
\begin{align}
|\varphi^n(AC)| &= \#\big((\Lambda^n t_k,\Lambda^n t_{k+2}]
  \cap \Sigma_{\varepsilon,\eta}(\Omega)\big)\,, \label{eq:vaphiAC}\\
|\varphi^n(B)| &= \#\big((\Lambda^n t_m,\Lambda^n t_{m+1}]
  \cap \Sigma_{\varepsilon,\eta}(\Omega)\big)\,. \label{eq:varphiB}
\end{align}

By definition of the matrix of a morphism and by
Corollary~\ref{coro:sumy-radku}, we have
$|\varphi^n(AC)|-|\varphi^n(B)|=\pm(\det\mat{M})^n$. Observe that
intervals $(\Lambda^n t_k,\Lambda^n t_{k+2}]$ and $(\Lambda^n
t_m,\Lambda^n t_{m+1}]$ have the same length, namely
$\Lambda^n(1+2\eta)$, and that the equality
$\lambda'\Z[\varepsilon']=\Z[\varepsilon']$ holds due
to~(\ref{eq:MF-ML}). We can therefore use Lemma~\ref{lem:Pk},
which states that the difference between the right hand sides
of~(\ref{eq:vaphiAC}) and~(\ref{eq:varphiB}) is bounded by a
constant $R$ independent of $n$.  Putting both facts together one
obtains
\[
|\det\mat{M}^n| \leq R\quad\text{ for any $n\in\N$.}
\]
The statement follows from the fact that $\det\mat{M}$ is an integer.
\end{proof}

\begin{coro}\label{coro:primitive-staci}
The incidence matrix of a 3iet preserving morphism satisfies $|\det\mat{M}|\leq 1$.
\end{coro}
\begin{proof}
Consider the primitive morphism $\varphi_0$, defined in
Example~\ref{ex:primitive-3iet-pres}, and let us denote by $\mat{M}_0$ a
power of this matrix, which is positive. Let $\varphi$ be a
non-primitive 3iet preserving morphism and let $\mat{M}$ be its
matrix. The matrix $\mat{M}\mat{M}_0$ is positive, and thus it is
the matrix of a primitive 3iet preserving morphism. By
Proposition~\ref{prop:detM-01} we have
\[
1 \geq |\det(\mat{M}\mat{M}_0)| =
|\det\mat{M}|\underbrace{|\det\mat{M}_0|}_{=1} = |\det\mat{M}|\,.
\qedhere
\]
\end{proof}

\begin{thmb}
  Let $\varphi$ be a 3iet preserving morphism and let $\mat{M}$ be its incidence matrix.
  Then one of the following holds \vspace{-0.5em}
  \begin{itemize}
  \item
    $\det\mat{M}=0$ and $\varphi(u)$ is degenerated for every 3iet word $u$, \vspace{-0.5em}
  \item
    $\det\mat{M}=\pm 1$ and $\varphi(u)$ is non-degenerated for every non-degenerated 3iet
    word $u$.
  \end{itemize}
\end{thmb}
\begin{proof}
We use notation from Lemma~\ref{lem:podprostorP} and
Remark~\ref{ex:regular}. In particular, recall
\[
\vec{x}_1\coloneq
\Bigl(\begin{smallmatrix}1\\1\\0\end{smallmatrix}\Bigr)\,,\qquad
\vec{x}_2\coloneq
\Bigl(\begin{smallmatrix}0\\1\\1\end{smallmatrix}\Bigr)\,,\qquad
\vec{x}_3\coloneq
\Bigl(\begin{smallmatrix}0\\1\\0\end{smallmatrix}\Bigr)\,.
\]

Lemma~\ref{lem:podprostorP} states that $\mat{M}{\cal P}={\cal
P}$, where ${\cal P}$ is a subspace of $\R^3$ spanned by vectors
$\vec{x}_1$, $\vec{x}_2$. Let us denote by $\mathcal{S}$ the
lattice $\mathcal{S} = \Z\vec{x}_1 + \Z\vec{x}_2$. The action of
the matrix $\mat{M}$ on the 2-dimensional subspace $\mathcal{P}$
has the matrix $\widetilde{M}$ from~\eqref{eq:msvlnkou}, which has
by Lemma~\eqref{lem:e1e2e3} determinant $\delta\in\{1,-1\}$.
Therefore the vectors $\mat{M}\vec{x}_1$ and $\mat{M}\vec{x}_2$
form a basis of $\mathcal{S}$ as well and thus
\begin{equation}\label{eq:SMS}
\mat{M}\mathcal{S} = \mathcal{S}\,.
\end{equation}

By an easy computation, $\mat{M}\vec{x}_3 = m_{12}\vec{x}_1 +
m_{32}\vec{x}_2 +
  (-m_{12}-m_{32}+m_{22})\vec{x}_3$,
hence by~\eqref{eq:det_PMP} we have $\mat{M}\vec{x}_3\in\delta\Delta\vec{x}_3 +
\mathcal{S}$, where $\Delta=\det\mat{M}$ as before.
Moreover,
\[
\Bigl(\begin{smallmatrix}1\\1\\1\end{smallmatrix}\Bigr) =
   -\vec{x}_3 + \vec{x}_1 + \vec{x}_2 \quad\implies\quad
   \mat{M}\Bigl(\begin{smallmatrix}1\\1\\1\end{smallmatrix}\Bigr) \in
   -\delta\Delta\vec{x}_3 + \mathcal{S}\,,
\]
and if we replace $\vec{x}_3$ on the right-hand side using
$-\vec{x}_3 = \Bigl(\begin{smallmatrix}1\\1\\1\end{smallmatrix}\Bigr)-\vec{x}_1 -\vec{x}_2$
we obtain
\begin{equation}\label{eq:M1inS}
\mat{M}\Bigl(\begin{smallmatrix}1\\1\\1\end{smallmatrix}\Bigr) \in
   \delta\Delta\Bigl(\begin{smallmatrix}1\\1\\1\end{smallmatrix}\Bigr) + \mathcal{S}\,.
\end{equation}

\underline{Case 1:} Let $\Delta = \det\mat{M} = 0$. Then
by~\eqref{eq:M1inS} and~(\ref{eq:SMS}) there exist $K_1,L_1\in\Z$
such that
\[
\mat{M}\Bigl(\!\begin{smallmatrix}1\\1\\1\end{smallmatrix}\Bigr) =
\mat{M}\Bigl(\!\!\begin{smallmatrix}K_1\\K_1+L_1\\L_1\end{smallmatrix}\!\!\Bigr)
\quad\hbox{ which implies }\quad
(\alpha,\beta,\gamma)\mat{M}\Bigl(\!\begin{smallmatrix}1\\1\\1\end{smallmatrix}\Bigr)
= (\alpha,\beta,\gamma)
\mat{M}\Bigl(\!\!\begin{smallmatrix}K_1\\K_1+L_1\\L_1\end{smallmatrix}\!\!\Bigr)\,,
\]
for arbitrary parameters $(\alpha,\beta,\gamma)$. It means that
$(\alpha,\beta,\gamma)\mat{M}$ are parameters of a degenerated 3iet word.

\medskip

\underline{Case 2:} Let $\Delta = \det\mat{M} = \pm1$. Again,
by~\eqref{eq:M1inS} there exist $K_2,L_2\in\Z$ such that
\begin{equation}\label{eq:det1KL}
\mat{M}\Bigl(\!\begin{smallmatrix}1\\1\\1\end{smallmatrix}\Bigr) =
\pm\Bigl(\!\begin{smallmatrix}1\\1\\1\end{smallmatrix}\Bigr)
+ \Bigl(\!\!\begin{smallmatrix}K_2\\K_2+L_2\\L_2\end{smallmatrix}\!\!\Bigr)\,.
\end{equation}

We show that parameters $(\alpha,\beta,\gamma)\mat{M}$ correspond to a degenerated
3iet word only if the original parameters $(\alpha,\beta,\gamma)$ correspond to
a degenerated 3iet word.

Let $(\alpha,\beta,\gamma)$ be such that $(\alpha,\beta,\gamma)\mat{M}$ are
parameters of a degenerated 3iet word, i.e., there exist $K_3,L_3,H,S\in\Z$
\begin{equation}\label{eq:thm2-pstrany1}
(\alpha,\beta,\gamma)\mat{M}\Bigl(\!\begin{smallmatrix}1\\1\\1\end{smallmatrix}\Bigr)
= (\alpha,\beta,\gamma)
\mat{M}\Bigl(\!\!\begin{smallmatrix}K_3\\K_3+L_3\\L_3\end{smallmatrix}\!\!\Bigr) =
(\alpha,\beta,\gamma)
\Bigl(\!\!\begin{smallmatrix}H\\H+S\\S\end{smallmatrix}\!\!\Bigr)\,,
\end{equation}
where the last equality comes from~\eqref{eq:SMS}.
Multiplying the equation~\eqref{eq:det1KL} by $(\alpha,\beta,\gamma)$ from the left
one obtains
\begin{equation}\label{eq:thm2-pstrany2}
(\alpha,\beta,\gamma)\mat{M}\Bigl(\!\begin{smallmatrix}1\\1\\1\end{smallmatrix}\Bigr) =
\pm(\alpha,\beta,\gamma)\Bigl(\!\begin{smallmatrix}1\\1\\1\end{smallmatrix}\Bigr)
+ (\alpha,\beta,\gamma)
\Bigl(\!\!\begin{smallmatrix}K_2\\K_2+L_2\\L_2\end{smallmatrix}\!\!\Bigr)\,.
\end{equation}
Finally, comparing right-hand sides of~\eqref{eq:thm2-pstrany1}
and~\eqref{eq:thm2-pstrany2} we have
\[
(\alpha,\beta,\gamma)\Bigl(\!\begin{smallmatrix}1\\1\\1\end{smallmatrix}\Bigr) =
\pm(\alpha,\beta,\gamma)
\biggl(\!\!\begin{smallmatrix}K_2-H\\K_2+L_2-H-S\\L_2-S\end{smallmatrix}\!\!\biggr)\,,
\]
which means that
$(\alpha,\beta,\gamma)$ are parameters of a degenerated 3iet word.
\end{proof}


\section{Comments and open problems}

\begin{enumerate}[1)]
\item
  We have derived that matrices of 3iet preserving morphisms belong to the monoid
  $\text{E}(3,\N) := \{\boldsymbol{M}\in\mathbb{N}^{3\times 3}\;|\;
  \boldsymbol{M}\boldsymbol{E}\boldsymbol{M}^T = \pm\boldsymbol{E}\:
  \text{ and }\det\boldsymbol{M}=\pm 1\}$ where $\boldsymbol{E} =
  \Big(\!\begin{smallmatrix}0&1&1\\-1&0&1\\-1&-1&0\end{smallmatrix}\Big)$.
  Unfortunately --- in contrast to the sturmian case --- the opposite is not
  true. In fact, the monoid $\text{E}(3,\N)$ contains matrices associated with morphisms,
  which are not 3iet preserving. As an example one can consider the matrix
  $\mat{M} =\Big(\begin{smallmatrix}0&2&1\\2&3&5\\3&0&5\end{smallmatrix}\Big)$.
\item
  The mapping $\varphi\rightarrow\mat{M}_{\varphi}$, where $\varphi$ a is
  3iet preserving morphism and $\mat{M}_{\varphi}$ is its incidence matrix is not
  one-to-one. One can show that for
  $\big(\begin{smallmatrix}a&b\\c&d\end{smallmatrix}\big)
  \in\N^{2\times2}$ with $ad-bc=\pm1$ there exist $a+b+c+d-1$ different
  sturmian morphisms. The same question for matrices of 3iet preserving
  morphisms is not solved.
\item
  Unlike the free monoid
  $\text{SL}(2,\N)=\{\mat{M}\in\N^{2\times2}\:|\:\det\mat{M}=1\}$, which is generated
  by two matrices $\big(\begin{smallmatrix}1&1\\0&1\end{smallmatrix}\big)$
  and $\big(\begin{smallmatrix}1&0\\1&1\end{smallmatrix}\big)$, the monoid
  $\text{SL}(3,\N)$ is not free, and, moreover, it is not finitely
  generated~\cite[Appendix~A]{fogg}. It would be interesting to derive similar
  results for the monoid $\text{E}(3,\N)$.
\item
  Even though the aim of this paper is not to investigate explicite prescriptions
  of 3iet preserving morphisms, we can still provide some
  information about it, based on our results.
  It follows from the proof of Lemma~\ref{lem:e1e2e3} that for every 3iet
  preserving morphism $\varphi:\{A,B,C\}^*\rightarrow\{A,B,C\}^*$ the morphism
  given by $A\mapsto\sigma_{A,B}\circ\varphi(A)$ and
  $B\mapsto\sigma_{A,B}\circ\varphi(B)$ is sturmian, where
  $\sigma_{A,B}$: $A\mapsto A$, $B\mapsto AB$, $C\mapsto B$; analogously
  for morphisms $\sigma_{A,C}$ and $\sigma_{B,C}$.
\item
  In this paper we were not at all interested in the characterization of 3iet words,
  which are fixed points of primitive morphisms, that is, of 3iet words $u$ such
  that there exists a primitive morphisms $\varphi$ for which $\varphi(u)=u$.
  This question is completely solved for sturmian
  words~\cite{yasutomi-dm-2,balazi-masakova-pelantova-i-5}.
  Adamczewski~\cite{adamczewsi-jtnb-14} studied for 3iet words a weaker
  property, the so-called primitive substitutivity. An infinite word $u$ over an alphabet
  $\Acal$ is said to be primitively substitutive if there exists a word $v$ over an
  alphabet $\mathcal{B}$, which is a fixed point of a primitive morphism, and
  a morphism $\psi:\mathcal{B}^*\rightarrow\Acal^*$ such that $\psi(v)=u$.
  Adamczewski, using results of Boshernitzan and Carroll~\cite{boshernitzan-carroll-jam-72},
  proved that a non-degenerated 3iet word is primitively substitutive
  if and only if normalized parameters $\varepsilon,l,c$ (see~\eqref{eq:prevod})
  of the corresponding transformation belong to the
  same quadratic field. Similar study can be foud in~\cite{harriss-lamb-arxiv06}.
\end{enumerate}

\section*{Acknowledgements}

The authors acknowledge financial support by Czech Science
Foundation GA \v{C}R 201/05/0169, and by the grant LC06002 of the
Ministry of Education, Youth, and Sports of the Czech Republic.


\appendix

\section{Proof of Lemma~\ref{lem:Pk}}

In this Appendix we prove Lemma~\ref{lem:Pk}, which is rather
technical. The proof uses the four following claims.

\begin{app-claim}
  Let $\varepsilon,\eta$ be irrational numbers, $\varepsilon\neq-\eta$,
  and let $\Omega_1$, $\Omega_2$
  be arbitrary bounded intervals. Then
  $\#\big(\Omega_1 \cap \Sigma_{\varepsilon,\eta}(\Omega_2)\big) =
    \#\big(\Omega_2 \cap \Sigma_{\eta,\varepsilon}(\Omega_1)\big)$.
\end{app-claim}
\begin{proof}
\begin{align*}
\#\big(\Omega_1 \cap \Sigma_{\varepsilon,\eta}(\Omega_2)\big) &=
  \#\{a+b\eta\:|\:a,b\in\Z,\ a+b\eta\in\Omega_1,\ a-b\varepsilon\in\Omega_2\} = \\
  &= \#\{a+c\varepsilon\:|\:a,c\in\Z,\ a+c\varepsilon\in\Omega_2,\ a-c\eta\in\Omega_1\} = \\
  &= \#\big(\Omega_2 \cap \Sigma_{\eta,\varepsilon}(\Omega_1)\big)\,.
\qedhere
\end{align*}
\end{proof}

\begin{app-claim}
Let $\varepsilon$ be a quadratic irrational number with conjugate
$\varepsilon'$. Let $\lambda$ be a quadratic unit whose conjugate
$\lambda'$ satisfies $\lambda'\Z[\varepsilon'] =
\Z[\varepsilon']$. Then
$\lambda'\Sigma_{\varepsilon,-\varepsilon'}(\Omega) =
  \Sigma_{\varepsilon,-\varepsilon'}(\lambda\Omega)$.
\end{app-claim}

\begin{proof}
Let us consider $x = a-b\varepsilon\in\Z[\varepsilon]$. If we
denote $\eta = -\varepsilon'$, the number
$a+b\eta=a-b\varepsilon'\in\Z[\varepsilon']$ is the image of $x$
under the Galois automorphism and therefore we denote it by
$x'=a-b\varepsilon'$.

Note that the condition $\lambda'\Z[\varepsilon'] = \Z[\varepsilon']$
is equivalent to the condition $\lambda\Z[\varepsilon] = \Z[\varepsilon]$,
and that these two equalities imply that the mappings $x'\mapsto\lambda'x'$ and
$x\mapsto\lambda x$ are bijections on $\Z[\varepsilon']$ and $\Z[\varepsilon]$,
respectively.

By definition of a C\&P set we have
\[
\Sigma_{\varepsilon,-\varepsilon'}(\Omega) =
  \{x'\in\Z[\varepsilon']\:|\:x\in\Omega\}\,.
\]
We derive
\begin{align*}
\lambda'\Sigma_{\varepsilon,-\varepsilon'}(\Omega) &=
  \lambda'\{x'\in\Z[\varepsilon']\:|\:x\in\Omega\} =
  \{\lambda'x'\in\Z[\varepsilon']\:|\:\lambda x\in\lambda\Omega\} = \\
  &= \{y'\in\Z[\varepsilon']\:|\:y\in\lambda\Omega\} =
  \Sigma_{\varepsilon,-\varepsilon'}(\lambda\Omega)\,.
\qedhere
\end{align*}
\end{proof}

The following Claim is given without the proof, since it is just a special
case of Proposition~6.2 in~\cite{gmp-jtnb-15}.

\begin{app-claim}
Let $\hat{\varepsilon},\hat{\eta}$ be irrational numbers,
$\hat{\varepsilon}\neq-\hat{\eta}$, and let $\hat{\Omega}$ be an
arbitrary bounded interval. Then
\[
\Sigma_{\hat{\varepsilon},\hat{\eta}}((1+2\hat{\varepsilon})\hat{\Omega})
=
  (1-2\hat{\eta})\Sigma_{\frac{\hat{\varepsilon}}{1-2\hat{\varepsilon}},
  \frac{\hat{\eta}}{1+2\hat{\eta}}}(\hat{\Omega})\,.
\]
\end{app-claim}

\begin{app-claim}
Let $\tilde{\varepsilon},\tilde{\eta}$ be irrational numbers such that
$\tilde{\varepsilon}\neq-\tilde{\eta}$. Let $z\in\R$ and let $J$ be a bounded
interval. We denote
$Q(J,z)\coloneq\#(J\cap\Sigma_{\tilde{\varepsilon},\tilde{\eta}}(z-1,z])$.
Then there is a constant $R$ such that $|Q(J,z)-Q(J,t)|\leq R$ for every
$z,t\in\R$ and for every interval $J$.
\end{app-claim}
\begin{proof}
The condition $a-b\tilde{\varepsilon}\in(z-1,z]$, where
$a,b\in\Z$, can be equivalently rewritten as $a=\lfloor
z+b\tilde{\varepsilon}\rfloor=z+b\tilde{\varepsilon}-\{z+b\tilde{\varepsilon}\}$.
Hence
\[
\Sigma_{\tilde{\varepsilon},\tilde{\eta}}(z-1,z] =
\{b(\tilde{\varepsilon}+\tilde{\eta})+z-\{z+b\tilde{\varepsilon}\}\:|\:b\in\Z\}\,.
\]

We consider the interval $J$ with boundary points $c,c+l$, where
$c,l\in\R$ and $l>0$. If the point
$b(\tilde{\varepsilon}+\tilde{\eta})+z-\{z+b\tilde{\varepsilon}\}$
belongs to the set
$J\cap\Sigma_{\tilde{\varepsilon},\tilde{\eta}}(z-1,z]$, then
$c-z\leq b(\tilde{\varepsilon}+\tilde{\eta})\leq c+l-z+1$. On the
other hand, if $c-z+1<b(\tilde{\varepsilon}+\tilde{\eta})<c+l-z$
then the point
$b(\tilde{\varepsilon}+\tilde{\eta})+z-\{z+b\tilde{\varepsilon}\}$
belongs to
$J\cap\Sigma_{\tilde{\varepsilon},\tilde{\eta}}(z-1,z]$. It means
that the number of points in the set
$J\cap\Sigma_{\tilde{\varepsilon},\tilde{\eta}}(z-1,z]$ is at
least
$\big\lfloor\frac{l-1}{\tilde{\varepsilon}+\tilde{\eta}}\big\rfloor$
and at most
$\big\lceil\frac{l+1}{\tilde{\varepsilon}+\tilde{\eta}}\big\rceil$,
and hence
\[
\bigg\lfloor\frac{l-1}{\tilde{\varepsilon}+\tilde{\eta}}\bigg\rfloor \leq
Q(J,z) \leq \bigg\lceil\frac{l+1}{\tilde{\varepsilon}+\tilde{\eta}}\bigg\rceil\,.
\]
Note that the bounds on $Q(J,z)$ do not depend on $z$, and thus
the same estimate holds for $Q(J,t)$. Therefore
\[
|Q(J,z)-Q(J,t)| \leq
\bigg\lceil\frac{l+1}{\tilde{\varepsilon}+\tilde{\eta}}\bigg\rceil
-
\bigg\lfloor\frac{l-1}{\tilde{\varepsilon}+\tilde{\eta}}\bigg\rfloor\leq
2\bigg(1+\frac{1}{\tilde{\varepsilon}+\tilde{\eta}}\bigg)\eqcolon
R\,.
\]
\end{proof}

Now we are in the position to conclude the proof of
Lemma~\ref{lem:Pk}.
\begin{proof}[Proof of Lemma~\ref{lem:Pk}]
Let us recall the definition of $P_n(x) =
\#\Big(\big(x,x+(1+2\eta)\Lambda^n\big]\,\cap\,\Sigma_{\varepsilon,\eta}(\Omega)\Big)$.
By Claim 1, we have
\[
P_n(x) = \#\Big(\Omega \,\cap\,
\Sigma_{\eta,\varepsilon}(x,x+(1+2\eta)\Lambda^n]\Big) =
\#\Big(\Lambda^n\Omega \,\cap\,
\Lambda^n\Sigma_{\eta,\varepsilon}(x,x+(1+2\eta)\Lambda^n]\Big)\,.
\]
As $\eta=-\varepsilon'$ and $\lambda'=\Lambda$, we have, by Claim
2,
\[
P_n(x) = \#\Big(\Lambda^n\Omega \,\cap\,
\Sigma_{\eta,\varepsilon}\big(\lambda^nx,\lambda^nx+(1+2\eta)\big]\Big)\,,
\]
where we used $\lambda\lambda'=\lambda\Lambda = 1$. Claim 3
further implies
\[
P_n(x) = \#\bigg(\Lambda^n\Omega \,\cap\,
(1-2\varepsilon)\Sigma_{\frac{\eta}{1-2\eta},\frac{\varepsilon}{1+2\varepsilon}}
\Big(\frac{\lambda^nx}{1+2\eta},\frac{\lambda^nx}{1+2\eta}+1\Big]\bigg)\,.
\]
Thus $P_n(x)=Q\big(J,\frac{\lambda^nx}{1+2\eta}+1\big)$ as defined
in Claim 4, where $J=\frac{1}{1-2\varepsilon}\Lambda^n\Omega$,
$\tilde{\varepsilon}= \frac{\eta}{1-2\eta}$ and
$\tilde{\eta}=\frac{\varepsilon}{1+2\varepsilon}$. The statement
of Lemma follows by application of Claim 4.
\end{proof}


\end{document}